\documentclass[11pt,letterpaper]{article}
\usepackage{mymacros}

\begin{document}

\begin{titlepage}

\renewcommand{\thefootnote}{\fnsymbol{footnote}}
\pagenumbering{gobble}
\thispagestyle{empty}

\vspace*{2cm}

\begin{center}
\textbf{\LARGE 
A Nonsmooth Dynamical Systems 
Perspective\\[.5em] on Accelerated Extensions of ADMM
}

\vspace{1.0cm}

\textbf{Guilherme~Fran\c ca,${}^{\!ab}$\footnote{\mailto{guifranca@gmail.com}}}
\textbf{~Daniel~P.~Robinson,${}^{\!c}$\footnote{\mailto{daniel.p.robinson@gmail.com}}}
\textbf{~and Ren\' e~Vidal${}^{b}$\footnote{\mailto{rvidal@jhu.edu}}}

\vspace{.1cm}

${}^{a}$\textit{University of California, Berkeley}\\[.2cm]
${}^{b}$\textit{Mathematical Institute for Data Science,
Johns Hopkins University}\\[.2cm]
${}^{c}$\textit{Lehigh University}

\end{center}

\vspace{0.5cm}

\begin{abstract}
Recently, there has been great interest in connections
between continuous-time dynamical systems and optimization methods, 
notably in the context of accelerated  methods for smooth and unconstrained
problems. In this paper we extend this perspective to nonsmooth and
constrained problems by obtaining differential
inclusions associated to novel accelerated variants of the  
alternating direction method of multipliers (ADMM).
Through a  Lyapunov analysis, we derive rates of convergence
for these dynamical systems in different settings that
illustrate an interesting tradeoff between decaying
versus constant damping strategies. 
We also obtain modified equations capturing fine-grained details of these methods,
which have improved stability and preserve the leading order convergence rates.
An extension to general nonlinear equality and inequality constraints in connection
with singular perturbation theory is provided.
\end{abstract}

\end{titlepage}

\renewcommand\contentsname{}                                                   
\tableofcontents

\renewcommand*{\thefootnote}{\arabic{footnote}}
\setcounter{footnote}{0}
\pagenumbering{arabic}
\setcounter{page}{1}

\section{Introduction}

The simplest method to minimize a smooth function
$\objective : \mathbb{R}^n \to \mathbb{R}$ is 
\emph{gradient descent}, given by
$x_{k+1} = x_k - \epsilon \nabla \objective(x_k)$,
where $k=0,1,\dotsc$ is the iteration time and $\epsilon > 0$ is the step size.  
When $\objective$ is convex, gradient descent  converges 
at a rate  of $\order(1/ \epsilon k)$  \cite{Nesterov:2004}, and
when $\objective$ is $\mu$-strongly convex
it converges at a rate of $\order(e^{-\mu \epsilon k})$
\cite{Nesterov:2004}.%
\footnote{By convergence rate we mean
an upper bound $f(x_k) - f^\star = \order(1/\epsilon k)$
for convex functions, where $f^\star \equiv \min_x \objective(x)$.
For strongly convex functions one can  bound the trajectories,
$\| x_k - x^\star\|^2 = \order(e^{-\mu \epsilon k})$, where
$x^\star \equiv \argmin_x \objective(x)$;
$\mu > 0$
measures the curvature of $\objective$
(see Eq. \eqref{eq:strongly_convex} below).}
Clearly, gradient descent is an explicit Euler discretization
of the \emph{gradient flow},
\begin{equation}
\label{eq:gradient_ode}
\dot{X}(t) = - \nabla \objective(X(t)) ,
\end{equation}
where 
$x_k \approx X(t_k)$,  $t_k = \epsilon k $, and
$\dot{X} \equiv dX / dt$.
One can show that system \eqref{eq:gradient_ode}
has a convergence rate of $\order(1/t)$ for convex  $\objective$,
and a rate of $\order(e^{-\mu t})$ for $\mu$-strongly convex $\objective$,
matching the behavior of its discrete counterpart---these results follow as particular cases of
Proposition~\ref{th:relaxed_admm_convex}.
Although gradient flow/descent is widely used in machine learning
(see, e.g., \cite{Fattahi:2019} for applications in online and nonconvex
settings) and controls (see, e.g., \cite{Casau:2020}),
its convergence can be slow.

A popular approach to speed up
gradient descent was proposed by Nesterov
\cite{Nesterov:1983}. Such a method 
converges at a rate of $\order(1/\epsilon k^2)$
for a convex function \cite{Nesterov:1983}, which is known to be
optimal in the sense of worst-case complexity \cite{Nesterov:2004};
generalizations of Nesterov's method to monotone inclusions also exist~\cite{Monteiro:2013}.
``Acceleration'' in the context of optimization
has  been considered counterintuitive, without a clear understanding of its underpinning mechanism. Recently, the continuous-time limit of 
Nesterov's method was obtained as~\cite{Candes:2016}
\begin{equation}
\label{eq:nesterov_ode}
\ddot{X}(t) + \dfrac{\rr}{t} \dot{X}(t) = - \nabla \objective(X(t))
\qquad (\rr = 3).
\end{equation}
This connection has been useful in providing insights into accelerated methods.
Followup work has brought a larger class
of  methods into a Hamiltonian formalism \cite{Wibisono:2016},
analyzes based on Lyapunov's theory have  been explored
\cite{Krichene:2015,Wilson:2016}, and
connections with symplectic geometry and geometric integrators
have been established \cite{Franca:2019,franca2020dissipative,Franca:2021}.
However, such a fruitful interplay 
has been limited mostly to gradient-based methods, i.e.,
for smooth and unconstrained problems.\footnote{However, the first-order
gradient flow for constrained problems has been considered long ago
\cite{Yamashita:1980,Schropp:2000}}

An important algorithm to solve nonsmooth and
composite problems is the
\emph{alternating direction method of multipliers} (ADMM)
\cite{Gabay:1976,Glowinsky:1975};
see  \cite{Boyd:2011} for a review.
In the convex case, ADMM converges at a rate of $\order(1/\epsilon k)$
\cite{He:2012,Eckstein:2015}, while in the strongly convex case it converges
exponentially \cite{Deng:2016}.
Variants of ADMM exist, including one that
uses relaxation \cite{Eckstein:1994,Eckstein:1998}.
Relaxed ADMM (R-ADMM) also
has exponential convergence for strongly convex
functions \cite{Giselsson:2014,FrancaBento:2016}.
Moreover, an \emph{accelerated} variant of ADMM (A-ADMM) has been proposed
\cite{Goldstein:2014}.
While numerical experiments \cite{Goldstein:2014} show that A-ADMM may
outperform Nesterov's method, its convergence rates in general settings are unknown.

Recently, we considered the continuous-time 
limit of ADMM and A-ADMM
\cite{Franca:2018} and a preliminary Lyapunov 
analysis. 
Here 
we extend this perspective significantly.
Focusing on problem\footnote{The more general formulation
$\min_{x,z}\left\{ f(x) + g(z)
\, | \, \A x + \B z = c\right\}$ is equivalent
if $\B$ is invertible; indeed, since $z - \B^{-1}c = - \B^{-1}\A x$,
one can redefine $\A$ and translate $z$ to obtain a form
similar to \eqref{eq:minimize}.}
\begin{equation}
\label{eq:minimize}
\min_{x}  \left\{ \objective(x) \equiv f(x) + g(\A x) \right\},
\end{equation}
where
$f:\Reals^n\to\Reals \cup \{+\infty\}$
and
$g:\Reals^m\to\Reals \cup \{+\infty\}$
are \emph{nonsmooth}, and  $\A \in\Reals^{m\times n}$ has
full column rank, 
our main contributions are summarized as follows:
\begin{itemize}
\item We introduce new extensions of ADMM that
combine \emph{relaxation and acceleration} under two types
of damping, namely decaying and constant. 
We call such methods
\emph{relaxed and accelerated ADMM} (R-A-ADMM) and
\emph{relaxed heavy ball ADMM} (R-HB-ADMM), respectively.
\item We derive \emph{differential inclusions} modeling these methods
to leading order and
provide a (nonsmooth) Lyapunov analysis yielding convergence
rates in the convex and strongly convex settings.
Our results highlight a tradeoff between the type of
damping versus degree of convexity, 
as shown in Table~\ref{table:convergence}.\footnote{%
$\alpha$ is the relaxation parameter,
$\mu$ is the strong convexity parameter,
$\sigma_1(\A)$ is the largest singular value of $\A$,
$\rr$ is the constant in the damping coefficient.}
In the convex case, the 2nd-order system related to
R-A-ADMM (decaying damping) achieves the optimal rate, 
however in the strongly convex
case it is the 2nd-order system related to
R-HB-ADMM (constant damping) that achieves the optimal rate.
\item We provide  \emph{backward-error analysis} for the discrete-time methods, i.e., we
show that 
the modified equations describing these methods to next to leading order
have improved stability and preserve these convergence rates.
\item We extend our framework to incorporate general nonlinear equality and inequality
\emph{constraints}, besides establishing 
connections with \emph{multiple scale analysis} and \emph{singular
perturbation theory}.
\end{itemize}

\begin{table}[t]
\centering
\renewcommand{\arraystretch}{1.4}
\begin{tabular}{@{}llll@{}}
\toprule[1pt]
\makecell[l]{\emph{continuum system} } &
\makecell[l]{\emph{convex}} &
\makecell[l]{\emph{strongly convex}} &
\emph{proof} \\
\midrule[.5pt]
R-ADMM & $t^{-1}$ & $\exp\Big( - \tfrac{\alpha \mu }{\sigma_1^2(\A) } t\Big)$
& Prop.~\ref{th:relaxed_admm_convex}\\
R-A-ADMM
& $t^{-2}$ & $t^{-2\rr/3}$ & Prop.~\ref{th:relaxed_aadmm_convex}\\
R-HB-ADMM &
$t^{-1}$ & $\exp\bigg(- \sqrt{ \tfrac{ \alpha \mu   }{\sigma_1^2(\A) }  } t\bigg)$
& Prop.~\ref{th:heavy_ball_admm_convex}\\
\bottomrule[1pt]
\end{tabular}
\caption{
\label{table:convergence}
Convergence rates
of the continuous-time dynamical systems associated to the variants of
ADMM proposed in this paper.
}
\end{table}

We note in passing that our followup paper \cite{Franca:2019b} explores the reverse
direction than considered here, namelly how different  methods---accelerated ADMM, Douglas-Rachford, proximal point, proximal gradient, Tseng splitting, etc.---arise as
different discretizations of the same system;
ADMM corresponds to a rebalanced splitting technique from numerical analysis. 
We also explored
the stochastic regime via  Langevin dynamics.
However, no convergence rates  nor constraints were 
considered in \cite{Franca:2019b}.
The variants of ADMM proposed in this paper were also not considered therein.

This paper is organized as follows.
Sec.~\ref{sec:preliminaries} contains background material.
Sec.~\ref{sec:variants} introduces our accelerated variants
of ADMM and their dynamical systems, whose convergence rates are stated in
Sec.~\ref{sec:convergence}.
Sec.~\ref{sec:perturbed} discusses a more refined continuous-time limit
and backward-error analyses for the convergence rates.
The proofs of our main results are provided in Sec.~\ref{sec:proofs}.
Sec.~\ref{sec:numerical} contains supporting numerical results.
Our concluding remarks are in 
Sec.~\ref{sec:discussion}.
Some omitted proofs and the extension to general nonlinear constraints, in connection with singular
perturbation theory,
are presented in the Appendix.

\section{Preliminaries}
\label{sec:preliminaries}

\subsection{Notation}
For $x, y \in  \Reals^n$ let $\| x \|^2=x^T x$ be the $L_2$  norm 
and $\langle x, y \rangle = x^T y$ be the inner product.
 The $L_1$ norm is denoted as $\| x \|_1$.
Given $\A \in \mathbb{R}^{m\times n}$ we denote its largest and
smallest singular values by $\sigma_1(\A)$ and $\sigma_n(\A)$, respectively,
and its condition number by $\kappa(\A) \equiv \sigma_1(\A) / \sigma_n(\A)$.
The nuclear norm is  $\| \A \|_* \equiv \sum_i \sigma_i(\A)$.
The minimum of a function $\objective$ is denoted by
$\objective^\star \equiv \min_x \objective(x)$, and a minimizer by
$x^\star \equiv \argmin_x \objective(x)$.
The  asymptotic notation
$f(t) = \order(g(t))$ means that there exists a constant $M > 0$ and
a number $t_0$ such that $|f(t)| \le M g(t)$ for all $t \ge t_0$, where
$g(t) > 0$.

\subsection{Subdifferentials, convexity, and strong convexity}

Consider  $f:\mathbb{R}^n \to \mathbb{R}\cup\{+\infty\}$
with effective domain $\dom f \equiv \{ x \, \vert \, f(x) < \infty \}$. Its
subdifferential at $x\in\dom f$ is defined as 
$\partial f(x) \equiv \left\{ v \in \mathbb{R}^n \, \vert \,
f(y) - f(x) \ge \langle v, y - x \rangle   \right\}$
for all $y \in \dom f$ \cite{RockafellarBook}.
The subdifferential set
$\partial f(x)$ is always closed and convex for $x \in \interior(\dom f)$,
and if $f$ is convex it is  also nonempty \cite{RockafellarBook}.

\begin{definition}
\label{def:strongly_convex}
The function $f : \Reals^n \to \Reals \cup \{+\infty\}$
is $\mu$-strongly convex
if and only if there exists a constant $\mu > 0$ such that
\begin{equation}
\label{eq:strongly_convex}
f(y) \ge f(x) +
\langle v , y - x \rangle +
(\mu / 2) \| y - x\|^2
\end{equation}
for all $x,y \in \dom f$ and all $v \in\partial f(x)$.
The function $f$ is said to be convex if this holds with $\mu = 0$.
\end{definition}

A function $f : \Reals^n \to \Reals \cup \{+\infty\}$ is \emph{closed} if
its epigraph $\{ (x, c) \in \Reals^n \times \Reals \, | \, f(x) \le c \}$ is
a closed set.  The following conditions are assumed throughout the paper.

\begin{assumption}
\label{th:assumption}
The functions
$f: \Reals^n \to \Reals\cup \{+\infty\}$ and
$g: \Reals^m \to \Reals\cup \{+\infty\}$
in problem \eqref{eq:minimize} are
closed and  convex, and $\objective \equiv f(x) + g(\A x)$
has a bounded global minimum. The matrix $\A$ has full column rank (so that $\AtA$ is
invertible).
\end{assumption}

Thus, $f$ and $g$ are proper, lower semicontinuous, and 
$0 \in \interior(\dom g - \A \dom f)$.
Moreover, $\partial \objective $ is upper semicontinuous on
$\interior(\dom \objective)$ and \cite{BorweinBook}
\begin{equation}\label{lemma}
\partial (f + g \circ \A) (x) = \partial f(x) + \A^T \partial g(\A x).
\end{equation}
Therefore, Assumption \ref{th:assumption} ensures that standard 
 subdifferential  calculus apply; the subdifferentials
$\partial f$, $\partial g$ and $\partial \objective$ are
\emph{maximal monotone operators} \cite{BorweinBook,CellinaBook2}.

\subsection{Differential inclusions and continuous-time limit}

Let $\mathcal{I} \equiv [0,T] \subset \mathbb{R}_+$ and
$F: \mathcal{I} \times \mathbb{R}^n
\rightrightarrows  \mathbb{R}^n$ be a \emph{multivalued map}.
Consider the generic differential inclusion \cite{Zeidler2,CellinaBook2,Deimling,Camlibel:2021}
\begin{equation}
\label{eq:inclusion}
\dot{X}(t) \in F(t, X(t)), \qquad X(0) = x_0.
\end{equation}
By a solution of \eqref{eq:inclusion}
we mean an absolutely continous function
$X: \mathcal{I} \to \mathbb{R}^n$ which satisfies the inclusion
for almost all $t\in\mathcal{I}$.

\begin{theorem}[see \cite{CellinaBook2,Deimling}]\label{thm_existence}
Let $\mathcal{D} \in \mathbb{R}^n$ be closed and $F : \mathcal{I} \times \mathcal{D} \rightrightarrows \mathbb{R}^n$
be an upper semicontinuous set-valued mapping, with non-empty, compact, and
convex values. Suppose  $\|F(t, x)\| \le c(1+ \| x\|)$ for almost
all  $t\in \mathcal{I}$
and some constant $c > 0$.  Then for every $x_0 \in \mathcal{D}$
there exists a solution $X(t)$
of the differential inclusion \eqref{eq:inclusion}.
Moreover, either $T \to \infty$ or $X(t)$ tends to the boundary
of $\mathcal{D}$ as $t \to T$.
\end{theorem}

The Assumption~\ref{th:assumption} implies that $F = \partial \objective$
is upper semicontinuous, nonempty, closed- and convex-valued.
\footnote{This is equivalent to say that $F$ is maximal monotone \cite{BorweinBook,CellinaBook2,Zeidler2}.
}
In this paper, we assume the growth condition of
Theorem.~\ref{thm_existence}.
This ensures local solutions on sufficiently small time intervals.
Global existence holds provided $F$ does not blowup in finite time.
Such a growth condition also  implies
$\| \dot{X}(t) \| \le c(1 + \| X \|)$ for almost all $t \in \mathcal{I}$,
hence Gronwall's inequality \cite{Gronwall:1919,Bellman:1943} yields
a constant $c' > 1$ such that $\max_{t\in \mathcal{I}} \| X(t)\| \le c' - 1$,
i.e., solutions of \eqref{eq:inclusion} are \emph{uniformly
bounded} \cite{Deimling}.
Furthermore, under the following Lipschitz condition, it is known
that a solution of  \eqref{eq:inclusion}
exists and is unique \cite{Dontchev:1992,Lempio:1995}.

\begin{definition} \label{lipschitz_def}
A multivalued map
$F : \mathcal{I} \times \mathbb{R}^n \rightrightarrows \mathbb{R}^n$
satisfies the (one-sided) Lipschitz condition with contant $L > 0$ if
\begin{equation} \label{lip}
\langle \xi_{t,y} - \xi_{t,x} , y - x \rangle \le L \| y - x\|^2
\end{equation}
for all $x, y \in \mathbb{R}^n$ and all $t \in \mathcal{I}$, where
$\xi_{t,y} \in F(t,y)$ and $\xi_{t,x} \in F(t,x)$.
\end{definition}

Thus, we also assume the following throughout the paper.

\begin{assumption} \label{assump2}
The subdifferential $\partial \objective$
obeys the growth condition of Theorem~\ref{thm_existence} and the Lipschitz
condition \eqref{lip}.
\end{assumption}

Under the above  assumptions, it is known \cite{Dontchev:1992,Lempio:1995,Beyn:2010,Rieger:2014} that explicit, implicit, and
semi-implicit Euler discretizations of  the differential inclusion \eqref{eq:inclusion} 
yield a  solution, i.e.,
such discretizations replace the original
differential inclusion over  $\mathcal{I} \equiv [0,T]$ 
by a sequence of \emph{discrete inclusions} on a grid $0=t_0 < t_1 < \dotsm < t_N =T$, 
where $t_k \equiv \epsilon k$, $k=0,1,\dotsc,N$, 
and step  size $\epsilon \equiv T / N$.
Thus, for given $\epsilon > 0$,  
solutions of the discrete inclusions yield a sequence
$\{ x_k \}_{k=0}^N$ that \emph{converges uniformly} to a solution $X(t)$
of  \eqref{eq:inclusion} as $\epsilon \to 0$
\cite{Dontchev:1992,Lempio:1995,Beyn:2010,Rieger:2014}.
It is also possible to estimate the error in this approximation;
 for such methods and under the above
conditions,
$\| X(t_k) - x_k \| = \order(\epsilon)$ 
\cite{Dontchev:1992,Lempio:1995,Beyn:2010,Rieger:2014}.
Throughout, we use the notation  $x_k \approx X(t_k)$, i.e., 
we assume there exists a continuous function $X(t)$
that approximates $x_k$ at time $t_k$ such that this approximation
is exact when  $\epsilon \to 0$.
Thus, under our regularity assumptions, we can use  standard finite difference 
approximations,
$\dot{X}(t_k) = \pm ( x_{k\pm 1} - x_k ) / \epsilon +\bigO(\epsilon)$ 
and $\ddot{X}(t_k) = (x_{k+1} - 2x_k + x_{k-1}) / \epsilon^2 +\bigO(\epsilon)$,
when considering the continuous-time limit of a given discrete algorithm.

\section{Variants of ADMM and dynamical systems}
\label{sec:variants}

\subsection{Relaxed ADMM}

Define the \emph{proximal operator} of a function
$f : \mathbb{R}^n \to \mathbb{R}\cup\{+\infty\}$---%
 or equivalently the \emph{resolvent} of its subdifferential 
$\partial f$---by
\begin{equation} \label{eq:prox_def}
J_{\partial f}^{\A}(z) \equiv \argmin_{x \in \mathbb{R}^n} \left\{ 
f(x) + (1/2)\| \A x - z \|^2 \right\} ,
\end{equation}
where $\A \in \mathbb{R}^{m \times n}$.
If $\A = \bm{I}$ we denote this by 
$J_{\partial f}$.
Let us start by considering the known family of R-ADMM algorithms \cite{Boyd:2011}:
\begin{subequations}
\label{eq:relaxed_admm}
\begin{align}
x_{k+1} &\leftarrow
J_{\rho^{-1} \partial f}^{\A} \big( z_k - u_k\big), 
\label{eq:relaxed_admm1}
\\
z_{k+1} &\leftarrow
J_{\rho^{-1}\partial g}\big( \alpha \A x_{k+1} + (1-\alpha) z_k + u_k \big),
\label{eq:relaxed_admm2}
\\
u_{k+1} &\leftarrow u_k + \alpha \A x_{k+1} + (1-\alpha)z_k - z_{k+1}.
\label{eq:relaxed_admm3}
\end{align}
\end{subequations}
The relaxation parameter $\alpha
\in(1,2)$ may  improve convergence
\cite{Eckstein:1994,Eckstein:1998}. The standard ADMM method is
recovered with $\alpha=1$. Above, $\rho \equiv \epsilon^{-1} > 0$ 
is the penalty
parameter, or inverse step size.
The proof of the following result is presented in  Appendix~\ref{sec:der_inclusions}.

\begin{proposition}
\label{th:relaxed_admm_ode}
Consider algorithm \eqref{eq:relaxed_admm} for 
problem \eqref{eq:minimize} under Assumptions~\ref{th:assumption} and
\ref{assump2}.
Let $t \equiv \epsilon k$ ($k=0,1,\dotsc$). In the 
limit $\epsilon \equiv \rho^{-1} \to 0$ these
updates tend to 
\begin{equation}
\label{eq:relaxed_admm_ode}
\alpha^{-1} \AtA \dot X(t) \in
- \partial\objective(X(t)) ,
\end{equation}
where $X(0) = x_0 \in\dom\objective$ is the initial state.
\end{proposition}

We make a couple of remarks:

\begin{itemize}

\item When $f$ and $g$ are differentiable,
$\A = \bm{I}$, and $\alpha=1$, we recover the gradient flow
\eqref{eq:gradient_ode} from the differential inclusion \eqref{eq:relaxed_admm_ode}.
However, the term $\AtA$ can make
the stability of the system different, even in the smooth case.

\item The variable $u_k$ has no counterpart in \eqref{eq:relaxed_admm_ode} and
the functions $f$ and $g$ appear together through $\objective$.
Thus,  the splitting of $f$ and $g$ and the introduction
of $u_k$ correspond to a discretization technique. 
Indeed, in  \cite{Franca:2019b} it was shown that $u_k$ is  a
``balance coefficient'' whose role is
to preserve critical points of the system \eqref{eq:relaxed_admm_ode} when 
the discretization splits $f$ from $g$.
This is a complementary perspective to
the augmented Lagrangian approach, where $u_k$ appears as  a Lagrange multiplier~\cite{Boyd:2011}.

\end{itemize}

\subsection{Relaxed and accelerated ADMM}

We now introduce
variables $\hat{u}\in\Reals^m$ and $\hat{z}\in\Reals^m$ to propose
an accelerated version of algorithm \eqref{eq:relaxed_admm},
which we call R-A-ADMM with updates given by 
\begin{subequations}
\label{eq:relaxed_aadmm}
\begin{align}
x_{k+1} &\leftarrow
J_{\rho^{-1}\partial f}^{\A}\big(  \hat{z}_k - \hat{u}_k \big) ,
\label{eq:raadmm1} \\
z_{k+1} &\leftarrow 
J_{\rho^{-1}\partial g}\big( \alpha\A x_{k+1} + (1-\alpha)\hat{z}_k + \hat{u}_k \big), 
 \label{eq:raadmm2} \\
u_{k+1} &\leftarrow \hat{u}_k +
                \alpha \A x_{k+1} + (1-\alpha) \hat{z}_k - z_{k+1},
                \label{eq:raadmm3} \\
  \gamma_{k+1} &\leftarrow \big(k+\sqrt{\rho}\big)\big/\big(k+\rr+\sqrt{\rho}\big), \label{eq:raadmm33}\\
  \hat{u}_{k+1} &\leftarrow u_{k+1} +
                \gamma_{k+1} \left( u_{k+1}-u_k \right), \label{eq:raadmm4}\\
  \hat{z}_{k+1} &\leftarrow z_{k+1} +
                \gamma_{k+1} \left( z_{k+1}-z_k \right). \label{eq:raadmm5}
\end{align}
\end{subequations}
Compared to the accelerated method of\cite{Monteiro:2013},
the above method splits $f$ from $g$,  allowing both functions to be nonsmooth.
The following result (derived in Appendix~\ref{sec:der_inclusions})
shows that the associated continuous-time system to the above method is a generalization of \eqref{eq:nesterov_ode} to the nonsmooth and
linear constrained case (and does not suffer from a divergence when $t \to 0$).

\begin{proposition}
\label{th:relaxed_aadmm_ode}
Consider algorithm \eqref{eq:relaxed_aadmm} for 
problem \eqref{eq:minimize} under Assumptions~\ref{th:assumption} and
\ref{assump2}. Let $t \equiv \epsilon k$ ($k=0,1,\dotsc$). In the
limit $\epsilon \equiv \rho^{-1/2} \to 0$
such updates tend to 
\begin{equation}
\label{eq:relaxed_aadmm_ode}
\alpha^{-1} \AtA \left[ \ddot{X}(t) + \dfrac{\rr}{t+1} \dot{X}(t) \right]
\in  - \partial \objective(X(t))
\end{equation}
with initial conditions $X(0) = x_0 \in \dom\objective$ and $\dot{X}(0) = 0$.
\end{proposition}



\subsection{Relaxed heavy ball ADMM}

Another acceleration scheme for gradient descent is the heavy ball method
\cite{Polyak:1964}.  Motivated by this  we
introduce another accelerated variant  of ADMM, 
called R-HB-ADMM, with updates given by
\begin{subequations}
\label{eq:heavy_ball_admm}
\begin{align}
x_{k+1} &\leftarrow
J_{\rho^{-1}\partial f}^{\A}\big(  \hat{z}_k - \hat{u}_k \big) ,
\label{eq:heavy_ball_admm1}\\
z_{k+1} &\leftarrow 
J_{\rho^{-1}\partial g}\big( \alpha\A x_{k+1} + (1-\alpha)\hat{z}_k + \hat{u}_k \big), 
\label{eq:heavy_ball_admm2} \\
u_{k+1} &\leftarrow \hat{u}_k + \alpha \A x_{k+1} +
    (1-\alpha) \hat{z}_k - z_{k+1},\label{eq:heavy_ball_admm3}\\
\hat{u}_{k+1} &\leftarrow u_{k+1} + \gamma \left( u_{k+1}-u_k \right),
    \label{eq:heavy_ball_admm4}\\
\hat{z}_{k+1} &\leftarrow z_{k+1} + \gamma \left( z_{k+1}-z_k \right).
    \label{eq:heavy_ball_admm5}
\end{align}
\end{subequations}
These updates are essentially the same as 
\eqref{eq:relaxed_aadmm},
but now $\gamma_k$ is the specific constant
$\gamma \equiv 1 - \rr / \sqrt{\rho}$.

\begin{proposition}
\label{th:heavy_ball_admm_ode}
Consider algorithm \eqref{eq:heavy_ball_admm} for solving
problem \eqref{eq:minimize} under Assumptions~\ref{th:assumption}~and~\ref{assump2}. Let $t \equiv \epsilon  k$ ($k=0,1,\dotsc$). 
In the limit $\epsilon \equiv \rho^{-1/2} \to 0$ it reduces to
\begin{equation}
\label{eq:relaxed_heavy_ball_aadmm_ode}
\alpha^{-1} \AtA \left[   \ddot{X}(t) + \rr \dot{X}(t) \right]
\in  - \partial \objective(X(t)) 
\end{equation}
with initial conditions
$X(0)= x_0 \in \dom\objective$ and
$\dot{X}(0) = 0$.
\end{proposition}

The key difference between systems  \eqref{eq:relaxed_heavy_ball_aadmm_ode}
and \eqref{eq:relaxed_aadmm_ode} is the damping coefficient, 
which leads to different stability properties (e.g., as
shown in Table~\ref{table:convergence}).
For  system \eqref{eq:relaxed_aadmm_ode}
the damping vanishes asymptotically, thus 
strong oscillations may arise for large $t$, precluding asymptotic stability. However, 
system \eqref{eq:relaxed_heavy_ball_aadmm_ode} has a constant damping, which introduces
more dissipation, thus long term oscillations can be better controlled.
With that in mind, one could consider a combination of such terms, e.g.,
\begin{equation}
\gamma_{k+1} = \dfrac{ k+\sqrt{\rho} }{ k+r_1 +\sqrt{\rho} } + \left(1- \dfrac{r_2 }{ \sqrt{\rho} }\right)
\end{equation}
in place of update \eqref{eq:raadmm33}, for constants $r_1,r_2 > 0$. In the continuous-time limit
such an algorithm is associated to the differential inclusion
\begin{equation} \label{eq:di_combined}
\alpha^{-1} \AtA \left[ \ddot{X}(t) + \right(\dfrac{r_1}{t + 1} + r_2\left) \dot{X}(t) \right]
\in  - \partial \objective(X(t)) .
\end{equation}
Choosing small $\rr_2 > 0$ may thus improve attractiveness properties 
 for large $t$ compared to system~\eqref{eq:relaxed_aadmm_ode}.
Next, we make a few important remarks.

\begin{remark}
\label{rem:hamiltonian}
Systems \eqref{eq:relaxed_aadmm_ode}, 
\eqref{eq:relaxed_heavy_ball_aadmm_ode}, and \eqref{eq:di_combined}  are \emph{nonsmooth} Hamiltonian systems
\cite{Rockafellar:1970,Rockafellar:1991,Rockafellar:1994,Ioffe:1997}, with Hamiltonian
\begin{equation}
\label{eq:ham1}
H(X, P; t) \equiv 
\dfrac{1}{2} e^{-\eta(t)} \big\langle P, \M^{-1}P \big\rangle
  + \alpha e^{\eta(t)} \objective(X),
\end{equation}
where $\M \equiv \AtA$ plays the role of a ``mass matrix'' and suitable $\eta$, which
is a damping function (e.g., $\rr \log (t+1)$ or $\rr t$). 
The relaxation parameter $\alpha$ appears as a ``coupling constant.''
Hamilton's equations,
$\dot{X} \in \nabla_P H$ and $\dot{P} \in - \partial_X H$, yield
$\alpha^{-1}\M\big(  \ddot{X} + \dot{\eta} \dot{X} \big) \in
-\partial \objective(X)$,
of which \eqref{eq:relaxed_aadmm_ode},
\eqref{eq:relaxed_heavy_ball_aadmm_ode},~and~\eqref{eq:di_combined} are particular
cases.
Hamiltonian systems have an intrinsic \emph{symplectic structure},
which has been recently explored in optimization  \cite{Franca:2019,franca2020dissipative,Franca:2021}.
%
%
\end{remark}

\begin{remark}\label{rmk_constr1}
Systems \eqref{eq:relaxed_admm_ode},  \eqref{eq:relaxed_aadmm_ode}
and \eqref{eq:relaxed_heavy_ball_aadmm_ode} capture the \emph{linear
constraint} of ADMM explicitly. However, it is possible to incorporate
more general constraints as a penalty term in the objective function, e.g., 
through $g$, which are then  solved
by the proximal operator \eqref{eq:relaxed_admm2}/\eqref{eq:raadmm2}/\eqref{eq:heavy_ball_admm2}. The advantage of ADMM, 
enabled by the linear constraint,  is that it splits $f$ from  $g$.
For instance, consider the linear inequality constrained problem
\begin{equation}
\min_x \big\{ f(x) \ \vert \ \B x \le c \big\} 
\end{equation}
for $\B \in \mathbb{R}^{m \times n}$ and $c \in \mathbb{R}^m$.
This can be converted into the form \eqref{eq:minimize} by
introducing a slack variable $z \ge -c$ so that the 
constraint becomes $\B x + z = 0$.  Setting $\A = -\B$, and introducing
the characteristic function $g = \chi_{+}$, where 
$\chi_+(z) \equiv 0 $ if $z \ge -c$ and 
$\chi_+(z) \equiv +\infty$ otherwise, we have the equivalent problem
\begin{equation}
\min_x \big\{ f(x) + \chi_{+}(-\B x) \big \} ,
\end{equation}
which is readily suitable to the previous methods.
For instance, when applying R-HB-ADMM \eqref{eq:heavy_ball_admm} to this problem, 
in light of Proposition~\ref{th:heavy_ball_admm_ode},
we can obtain useful insights by studying  the  nonsmooth system
\begin{equation}
\alpha^{-1}\B^T \! \B\left[ \ddot{X} + \rr \dot{X} \right] \in -\partial f(X) 
+ \B^T \partial \chi_{+}(-\B X) .
\end{equation}
 The same comment applies to the other ADMM variants and their
associated dynamical systems. The convergence rates we derive in the next section
apply to this problem.
\end{remark}

\begin{remark} \label{rmk_constr2}
To illustrate the generality of the previous connections,
consider the constrained problem
\begin{equation}
\min_{x \in \C} f(x) ,
\end{equation}
where $\C \subset \mathbb{R}^n$.  
Using the characteristic function $\chi_{\C}(z) \equiv 0 $ if $z \in \C$ and
$\chi_{\C}(z) \equiv +\infty$ if $z \notin \C$, 
and setting $\A = \bm{I}$ and $g = \chi_{\C}$,
the proximal operator \eqref{eq:relaxed_admm2}/\eqref{eq:raadmm2}/\eqref{eq:heavy_ball_admm2}
becomes the projection into $\C$.
The previous dynamical systems provide valuable insights into the behavior
of the algorithms on this problem, e.g., for R-HB-ADMM one would study
\begin{equation}
\alpha^{-1}\left[ \ddot{X} + \rr \dot{X} \right] \in - \partial f(X) - \partial \chi_{\C}(X) ,
\end{equation}
and similarly for the other ADMM variants.
In particular, if $\C$ is a convex set, the convergence rates obtained in the next section
apply to this case.
We provide a more general discussion about nonlinear constraints in 
Appendix~\ref{sec:constrained} as well as connections with multiple scale analysis and
singular perturbation theory.
\end{remark}

\section{Convergence rates}
\label{sec:convergence}

We  state several convergence rate results for the above differential inclusions.
The proofs of these results are presented in Section~\ref{sec:proofs}.

\begin{proposition}
\label{th:relaxed_admm_convex}
Consider the differential inclusion \eqref{eq:relaxed_admm_ode},
which models the R-ADMM \eqref{eq:relaxed_admm} to leading order.
If the objective function $\objective$ is convex then for almost all $t \in \mathcal{I}\subset \mathbb{R}_+$
it holds that
\begin{equation}
\label{eq:relaxed_admm_convex}
\objective(X(t)) - \objective^\star \le
 \dfrac{ \sigma_1^2(\A) \| x_0-x^\star\|^2 }{ 2 \alpha t }  .
\end{equation}
If $\objective$ is $\mu$-strongly convex then  for
almost all $t \in \mathcal{I}$ we have
\begin{equation}
\label{eq:relaxed_admm_strongly_convex}
\|X(t) - x^\star \|^2 \le \kappa^2(\A) \| x_0 - x^\star \|^2 \exp\left( - \dfrac{\mu \alpha  t}{ \sigma_1^2(\A)} \right)  .
\end{equation}
\end{proposition}

The convergence rate \eqref{eq:relaxed_admm_convex}
matches the $\order(1/\epsilon k)$ rate of standard ADMM
in the convex case \cite{Eckstein:2015,He:2012}.
We believe the  result 
\eqref{eq:relaxed_admm_convex}
for \emph{relaxed} ADMM
is new in the sense that this rate
is unknown in discrete time.
The exponential convergence rate \eqref{eq:relaxed_admm_strongly_convex}
agrees with the linear convergence of vanilla and relaxed ADMM
\cite{Hong:2017,FrancaBento:2016}.


\begin{proposition}
\label{th:relaxed_aadmm_convex}
Consider the differential inclusion \eqref{eq:relaxed_aadmm_ode} with $r \geq 3$, which
models the R-A-ADMM method \eqref{eq:relaxed_aadmm} to leading order.  If $\objective$ is convex then
for almost all $t \in \mathcal{I} \subset \mathbb{R}_+ $ we have
\begin{equation}
\label{eq:relaxed_aadmm_convex}
\objective(X(t)) - \objective^\star \le
\dfrac{C}{(t+1)^{2}}
\end{equation}
for a constant $ C > 0$. 
If $\objective$ is $\mu$-strongly convex
then there exists $t_0 > 0$ 
and a constant $C \ge 0$ 
such that  for almost all $t \ge \max\{ 0, t_0\} \in \mathcal{I}$
we have
\begin{equation}
\label{eq:relaxed_aadmm_strongly_convex}
\| X(t) - x^\star \|^2 \le  \dfrac{ 4  C }{ \alpha \mu   (t+1)^{2\rr/3} } .
\end{equation}
\end{proposition}

According to \eqref{eq:relaxed_aadmm_strongly_convex},  the 
system \eqref{eq:relaxed_aadmm_ode} does not attain exponential convergence;
its decaying damping does not provide enough dissipation,
in contrast with \eqref{eq:relaxed_admm_strongly_convex} and our next result for system~\eqref{eq:relaxed_heavy_ball_aadmm_ode}.

\begin{proposition}
\label{th:heavy_ball_admm_convex}
Consider the differential inclusion \eqref{eq:relaxed_heavy_ball_aadmm_ode},
which models the R-HB-ADMM method~\eqref{eq:heavy_ball_admm} to leading order.
If the objective function $\objective$ is convex then
for almost all $t \ge t_0 \equiv 1/\rr \in \mathcal{I} \subset \mathbb{R}_+$ we have
\begin{equation}
\label{eq:heavy_ball_admm_convex}
\objective(X(t)) - \objective^\star \le \dfrac{ C }{ \alpha t } 
\end{equation}
for some constant $C > 0$. 
If $\objective$ is $\mu$-strongly convex then for
almost all $t \in \mathcal{I}$ we have
\begin{equation}
\label{eq:heavy_ball_admm_strongly_convex}
\| X(t) - x^\star\|^2
\le  \dfrac{C}{  \mu} \exp \left( - \dfrac{ 2 \rr t}{3} \right) 
\end{equation}
for some constant $C > 0$, 
provided the damping coefficient obeys
\begin{equation} \label{eq:rbardef}
\rr  \le \bar{\rr}
\equiv \dfrac{ 3  \sqrt{\alpha\mu}  }{2  \sigma_1(\A) }  .
\end{equation}
\end{proposition}

We  note a few points regarding the results obtained thus far:

\begin{itemize}

\item The relaxation parameter $\alpha$ improves
convergence in some cases more than others, e.g., in the convex case
\eqref{eq:relaxed_aadmm_convex} there is only a linear improvement in
the constant $C$, but
in the strongly convex case
\eqref{eq:heavy_ball_admm_strongly_convex} the parameter $\alpha$  appears
under an exponential. 

\item The dynamical system modelling R-HB-ADMM attains exponential
convergence \eqref{eq:heavy_ball_admm_strongly_convex}.
In contrast, the system associated to R-A-ADMM achieves only a sublinear
rate \eqref{eq:relaxed_aadmm_strongly_convex}.
We observed numerically that
R-HB-ADMM outperforms
R-A-ADMM in most cases, even for nonstrongly convex functions.

\item Proposition~\ref{th:heavy_ball_admm_convex} applies to Polyak's heavy ball
method \cite{Polyak:1964} as a particular case---with $\A = \Id$, $\alpha=1$,
and $f$ smooth (the function $g=0$ is absent). In this case, a discrete-time
convergence rate of $\order(1/\epsilon k)$ was obtained
only in an average sense \cite{Ghadimi:2014}.
In the strongly convex case, exponential convergence
is known \cite{Ghadimi:2014}, but not in the
form \eqref{eq:heavy_ball_admm_strongly_convex}.
We believe the rates of Proposition~\ref{th:heavy_ball_admm_convex} 
have yet no discrete-time analog.


\item Comparing Propositions~\ref{th:relaxed_aadmm_convex} and
\ref{th:heavy_ball_admm_convex}, we see
 an interesting  tradeoff between  
the type of damping versus degree of convexity of the objective function.
Decaying damping achieves the optimal rate of $\order(1/t^2)$ in the
convex case,
but only a sublinear rate in the strongly convex case.
On the other hand, constant damping has a suboptimal rate 
in the convex case, but  the optimal linear rate of $\order\big(e^{-\sqrt{\alpha\mu} t / \sigma_1(\A)}\big)$
in the strongly
convex case. This indicates that when
$\objective$ has ``sufficient curvature,'' algorithms based on the latter (constant damping)
are preferable.

\end{itemize}

\section{Modified or perturbed differential equations}
\label{sec:perturbed}

The systems  \eqref{eq:relaxed_admm_ode}, \eqref{eq:relaxed_aadmm_ode},
and \eqref{eq:relaxed_heavy_ball_aadmm_ode} are able to capture the \emph{leading order} behavior
of their associated algorithms; the discrepancy  is 
$\| X(t_k) - x_k \| = \bigO(\epsilon)$. 
Thus, they  describe the main trend
but not fine-grained details of their associated algorithms, 
such as dependence on the step size and curvature of the objective function. 

When the objective function is smooth, a more refined approximation to Nesterov's method
was obtained \cite{Shi:2022}, as well as the modified equation
from backward error analysis \cite{Franca:2019}.
A key ingredient  is  the appearance of a Hessian-dependent damping,
implying an improved stability 
 \cite{Poveda:2020}.
We now show that similar perturbed dynamical systems, i.e., containing a spurious
Hessian-driven dissipation, are related to the previous
accelerated  ADMM methods. Thus, such a system models the associated algorithm
up to a discrepancy of 
$\| X(t_k) - x_k \| = \bigO(\epsilon^2)$,  providing an accurate description.
As common in numerical analysis, to obtain  perturbed equations
it is necessary to  assume sufficient smoothness.

\begin{assumption} \label{th:assump4}
Besides Assumption~\ref{th:assumption}, here  we also assume that
$f$ and $g$ are twice continuously  differentiable  and have Lipschitz continuous
gradients. For the  stability results, we also assume that $\objective$ is convex and
radially unbounded
($\objective(x) \to \infty$ as $\|x\| \to \infty$).
\end{assumption}

The following result is derived in the next section.

\begin{proposition}\label{th:mod_eq_aadmm}
Under the smoothness condition of Assumption~\ref{th:assump4}, 
the perturbed  differential equations modeling 
algorithms
\eqref{eq:relaxed_aadmm} and \eqref{eq:heavy_ball_admm}, which  are
accurate up to $\bigO(\epsilon^2)$ ($\epsilon = \rho^{-1/2}$), 
have the form
\begin{equation}\label{eq:mod_eq_aadmm3}
\dfrac{\AtA}{\alpha} \left[ \ddot{X} + a(t) \dot{X} \right] 
+ \epsilon \big( \nabla^2\objective \big) \dot{X}  
= - \big(1 + \epsilon  b(t) \big) \nabla \objective  ,
\end{equation}
where $X \equiv X(t)$ and  $\objective \equiv \objective(X(t))$. 
In the case of the R-A-ADMM method \eqref{eq:relaxed_aadmm} 
we have the time-dependent coefficients
\begin{equation} \label{eq:ab_dec}
a(t) = \dfrac{\rr}{t+1} - \dfrac{\epsilon \rr (r - 2) }{2(t+1)^2}, \qquad 
b(t)  = \dfrac{\rr}{2(t+1)},
\end{equation}
while for the  R-HB-ADMM method \eqref{eq:heavy_ball_admm}  we have
the constant coefficients
\begin{equation} \label{eq:ab_const}
a(t) = \rr + \dfrac{ \epsilon \rr^2 }{ 2 } , \qquad
b(t) = \dfrac{ \rr }{2}.
\end{equation}
\end{proposition}

When $\epsilon \to 0$ we recover \eqref{eq:relaxed_aadmm_ode} and \eqref{eq:relaxed_heavy_ball_aadmm_ode} from 
\eqref{eq:mod_eq_aadmm3}. 
The crucial ingredient in  the latter is the 
term $\epsilon (\nabla^2 \objective ) \dot{X}$
that introduces a \emph{spurious dissipation}, rendering the system 
more stable, i.e., the
methods~\eqref{eq:relaxed_aadmm}
and \eqref{eq:heavy_ball_admm}  have even better stability/attractiviness 
properties compared to their associated leading order differential inclusions. 
Indeed, consider 
\begin{equation}\label{eq:lyap_pert}
\E(X, \dot{X}, t) \equiv \dfrac{1}{2\alpha}\big\| \A \dot{X}\big\|^2 + 
\big(1+\epsilon b(t)\big)\big(\objective(X) - \objective^\star\big).
\end{equation}
Note that $\E > 0$ outside a critical point $(X, \dot{X}) = (x^\star, 0)$.
Under the evolution of system \eqref{eq:mod_eq_aadmm3} we have
\begin{equation}\label{eq:der_lyap_pert}
\dot{\E} = - a(t) \alpha^{-1} \big\| \A \dot{X}\big\|^2 - 
\epsilon \big\langle \dot{X}, \nabla^2 \objective \dot{X} \big\rangle
+ \dot{b}(t)\big( \objective(X) - \objective^\star \big) \le 0,
\end{equation}
since $\dot{b}(t) \le 0$ and $\nabla^2 \objective \succcurlyeq 0$. This implies
global uniform stability. The term $-\epsilon \big\langle \dot{X}, \nabla^2 \objective \dot{X}\big\rangle$ introduces extra contraction and is particularly
important when $a(t), \dot{b}(t) \to 0$ as $t \to \infty$ for the
case~\eqref{eq:ab_dec}; i.e., we would have
$\dot{\E} \to 0$ as $t\to\infty$ without  such a term, precluding asymptotic stability.

In fact, a stronger notion called \emph{uniform global asymptotic stability} (UGAS) \cite[Def. 3]{Poveda:2020}, which
combines uniform stability and uniform attractivity, was thoroughly studied for systems
in the form~\eqref{eq:mod_eq_aadmm3} \cite{Poveda:2020}.
Via the ``Lyapunov'' function \eqref{eq:der_lyap_pert} this analysis extends 
immediately to  \eqref{eq:mod_eq_aadmm3}, whence one concludes:

\begin{corollary}
Under  Assumption~\ref{th:assump4} we have the following results.
The autonomous system \eqref{eq:mod_eq_aadmm3}/\eqref{eq:ab_const} is 
UGAS (this is true even when $\epsilon \to 0$); this 
follows from \cite[Thm. 1]{Poveda:2020}.
The system \eqref{eq:mod_eq_aadmm3}/\eqref{eq:ab_dec} is UGAS (here we must have
$\epsilon > 0$);
this follows from \cite[Thm. 3]{Poveda:2020}. 
\end{corollary}

Thus, the perturbed system \eqref{eq:mod_eq_aadmm3} introduces 
a beneficial spurious dissipation that not only improves stability but, as stated
below,   also preserves the previous  convergence rates.

\begin{proposition}\label{th:preserve_rates_pert}
The convergence rates of Proposition \ref{th:relaxed_aadmm_convex} remain
valid for the perturbed or modified system \eqref{eq:mod_eq_aadmm3}/\eqref{eq:ab_dec}.
Similarly, the convergence rates of Proposition \ref{th:heavy_ball_admm_convex}
hold true for the perturbed system \eqref{eq:mod_eq_aadmm3}/\eqref{eq:ab_const}.
\end{proposition}

These results provide strong evidence that the accelerated ADMM variants  we introduced
achieve optimal rates of convergence. 
This is because the perturbed system captures a step size dependency and properties
of $\objective$. Note that higher order terms in an asymptotic expansion can only contribute with
even smaller effects which cannot spoil these results (for a suitable choice of step size).
Thus, the above results constitute \emph{backward error analysis} for the discrete-time
methods, and are valid for finite step sizes $\epsilon > 0$ (i.e., not only
when $\epsilon \to 0$ as before).
Our proof strategy for Proposition~\ref{th:preserve_rates_pert}  also extend immediately 
to other methods in the literature, beyond those
considered in this paper.

\section{Proofs of the main results}
\label{sec:proofs}

\subsection{Derivation of the differential inclusions}

The proofs of Propositions~\ref{th:relaxed_admm_ode},
\ref{th:relaxed_aadmm_ode}, and \ref{th:heavy_ball_admm_ode} are simpler than
that of Proposition~\ref{th:mod_eq_aadmm}, which will be
presented in this section, 
therefore they are deferred to  Appendix~\ref{sec:der_inclusions}.

\subsection{Derivation of the convergence rates}
\label{sec:der_rates}

Under Assumptions~\ref{th:assumption}
and \ref{assump2}, our differential inclusions have a unique solution.
Moreover, such a solution can be obtained as the limiting trajectory  
of a \emph{regularized differential equation}.
To justify this procedure, let 
us  first introduce the necessary concepts.

Consider the differential inclusion \eqref{eq:inclusion}
where $F = \partial \objective$.
It is useful
to consider the \emph{Moreau-Yosida
regularization} \cite{Zeidler2,CellinaBook2}.
The
\emph{Moreau envelope} of a function $\objective$ is defined by
\begin{equation} \label{eq:moreau}
\objective_\lambda(x)  \equiv \min_{y}
\left\{
\objective(y) +  \dfrac{1}{2\lambda} \| y - x \|^2
\right\} 
\end{equation}
for all $x\in \dom \objective$,
where $\lambda > 0$.  An important property is that 
$\objective^\star \le \objective_\lambda(x) \le \objective(x)$, i.e., the Moreau envelope
lower bounds the function, and moreover
$\objective$ and $\objective_\lambda$ have exactly the same minimizers,
in which case
$\objective_\lambda(x^\star) = \objective(x^\star) \equiv \objective^\star$.
It is also well-known that
$\nabla \objective_{\lambda}(x) =
\lambda^{-1}\left(x - J_{\lambda \partial \objective}(x)\right)$,
where $J_{\lambda \partial \objective}$ is the proximal
operator defined in \eqref{eq:prox_def}.
When $\partial \objective$ is maximal monotone,
$J_{\lambda \partial \objective}$ 
is single-valued and $\nabla \objective_\lambda$ is Lipschitz continuous.
Thus, instead of the differential inclusion
$\dot{X}(t) \in -\partial \objective(X(t))$,
one can study the regularized differential equation
$\dot{X}_\lambda(t) = -\nabla \objective_\lambda(X_\lambda(t))$,
which is well-posed. Under our 
 assumptions,
$\lim_{\lambda \downarrow 0}X_\lambda(t) = X(t)$, where $X(t)$ is a solution
of the the differential inclusion \cite{Zeidler2,CellinaBook2}.

We can therefore replace 
\eqref{eq:relaxed_admm_ode}, \eqref{eq:relaxed_aadmm_ode} and
\eqref{eq:relaxed_heavy_ball_aadmm_ode} by their
regularized differential equations
\begin{align}
\alpha^{-1} \AtA \dot{X}_{\lambda}(t)
&= - \nabla \objective_\lambda(X_{\lambda}(t)), \label{ode1} \\
\alpha^{-1}\AtA \left[\ddot{X}_\lambda(t) \! + \!
\dfrac{\rr}{t+1} \dot{X}_\lambda(t) \right] &= -
\nabla \objective_\lambda(X_\lambda(t)), \label{ode2} \\
\alpha^{-1} \AtA \left[\ddot{X}_\lambda(t) +
\rr \dot{X}_\lambda(t) \right] &= -
\nabla \objective_\lambda(X_\lambda(t)) , \label{ode3}
\end{align}
respectively.
These initial value problems are well-posed since $\nabla \objective_\lambda$ is Lipschitz
continuous and the differential operators on the left-hand-side are bounded.  
One can show that
$\nabla \objective_\lambda(x) \to \partial \objective (x)$
as $\lambda \downarrow 0$ \cite{Zeidler2,CellinaBook2}---this limit
actually yields the element of minimal norm in the
set $\partial \objective(x)$---and
$X(t) \equiv \lim_{\lambda \downarrow 0} X_\lambda(t)$ is
a solution of the associated differential inclusion.
In our case, Assumptions~\ref{th:assumption} and \ref{assump2}
ensure that the differential inclusions have a unique solution.
In what follows, we  derive convergence
rates for the above \emph{differential equations}, then  
take the limit $\lambda \downarrow 0$ which implies 
the same result for the associated differential inclusion.

To simplify notation, we will often omit the parameter $\lambda$ and 
restore it when necessary. We will also make extensive use of the ``perturbed variables''
\begin{equation}\label{eq:perturb}
\tilde{X} \equiv X - x^\star, \qquad
\tilde{\objective} \equiv \objective - \objective^\star ,
\end{equation}
where it is implied that $X \equiv X(t)$, $\objective \equiv \objective(X(t))$, and so on.

\begin{proof}[Proof of Proposition~\ref{th:relaxed_admm_convex}]
The convex case consists of a minor modification of a result 
from \cite{Franca:2018}, thus it is presented in Appendix~\ref{sec:more_rates}.
For the strongly convex case,  note that 
strong convexity of $\objective$  (see Eq.~\eqref{eq:strongly_convex})  implies
\begin{equation}
\label{eq:strongly_convex_minimum}
\big\langle \nabla \tilde\objective , \tilde{X}  \big\rangle \ge
\dfrac{\mu}{ 2} \big\| \tilde{X}  \big\|^2 
\ge \dfrac{\mu}{2 \sigma_1^{2}(\A)}  \big\| \A \tilde{X}  \big\|^2 .
\end{equation}
Consider
\begin{equation}
\label{eq:lyapunov2_2}
\E(X)  \equiv \dfrac{1}{2} \big\| \A  \tilde{X} \big\|^2.
\end{equation}
Taking its total time derivative along trajectories of \eqref{ode1},
and using \eqref{eq:strongly_convex_minimum}, yield
\begin{equation}
\begin{split}
\dot{\E} 
\le - \dfrac{ \mu \alpha}{ \sigma_1^{2}(\A) } \E .
\end{split}
\end{equation}
From Gr\" onwall's inequality \cite{Gronwall:1919,Bellman:1943} we thus
have
$\E(X_\lambda(t)) \le \E(x_0) e^{- \eta t}$, with $\eta \equiv \mu \alpha / \sigma_1^2(\A)$.
From standard norm inequalities, and taking the limit $\lambda \downarrow 0$,
we obtain the upper bound \eqref{eq:relaxed_admm_strongly_convex}
for the  differential inclusion \eqref{eq:relaxed_admm_ode}.
\end{proof}

\begin{proof}[Proof of Proposition~\ref{th:relaxed_aadmm_convex}]
The convex cases is similar to the derivation in \cite{Franca:2018}, thus
it is presented in Appendix~\ref{sec:more_rates}.
The proof for the strongly convex case is more involved and
takes several steps.
First, consider
\begin{equation}
\label{eq:lyapunov4_2}
\E \equiv
\alpha (t+1)^{w}  \tilde{\objective}
 + \dfrac{1}{2}(t+1)^{w-2}\big\| \A\big( w \tilde{X}  + (t+1) \dot{\tilde{X}} \big)  \big\|^2
\end{equation}
with $w \equiv 2\rr/3 \ge 2$ (since $\rr \ge 3$).
After simplifications, its total time derivative along trajectories
of system \eqref{ode2} becomes
\begin{equation}
\begin{split}
\label{eq:dotealmost}
\dot{\E} &=
\alpha w (t+1)^{w-1} \tilde{\objective}   \\
& + w^2(w/2-1) (t+1)^{w - 3}  \big\| \A \tilde{X} \big\|^2 \\
&+ w \left( w/2 -1\right)
(t+1)^{w-2} \big\langle \tilde{X}, \AtA \dot{\tilde{X}} \big\rangle \\
& - \alpha w (t+1)^{w - 1}  \big\langle
\tilde{X} , \nabla \tilde\objective  \big\rangle .
\end{split}
\end{equation}
From strong convexity of $\objective$, i.e., relation \eqref{eq:strongly_convex}, 
we have
\begin{equation}
\label{eq:stronga}
\big\langle \nabla \tilde \objective , \tilde{X} \big\rangle  \ge \tilde{\objective}  +
 \dfrac{ \mu }{ 2 \sigma_1^{2}(\A) }    \big\| \A \tilde{X} \big\|^2 .
\end{equation}
Thus, $\dot{\E}$ is less or equal than
\begin{multline}
\label{eq:dote}
w (t+1)^{w-3} \left[
w\left(\dfrac{w}{2}-1\right) - \dfrac{1}{2} \alpha \mu (t+1)^2 \sigma_1^{-2}(\A)
\right] \big\| \A \tilde{X} \big\|^2 \\
+ w \left(\dfrac{w}{2}-1\right)
(t+1)^{w-2} \big\langle \tilde{X} ,
\AtA \dot{\tilde{X}} \big\rangle .
\end{multline}
Requiring the first term above  to be nonpositive demands
\begin{equation} \label{t0choice}
t \ge t_0, \qquad t_0 \equiv  \sigma_1(\A) 
\sqrt{  w(w-2) (\alpha\mu)^{-1}  } - 1 . 
\end{equation}
Under such condition, we can neglect the first term above to obtain
\begin{equation}
\label{eq:dotefinal}
\dot{\E} \le  \dfrac{w}{2}
\left( \dfrac{w}{2}-1\right) (t+1)^{w-2}
\dfrac{d}{dt} \big\| \A \tilde{X}  \big\|^2 .
\end{equation}
Using integration by parts,
\begin{multline}
\label{eq:integral_strong}
\E(X(t), \dot{X}(t), t)  - \E(X(t_0), \dot{X}(t_0), t_0)  + 
\dfrac{w}{2}\left(\dfrac{w}{2}-1\right)
\bigg\{
(t_0+1)^{w-2} \big\| \A  \tilde{X}(t_0) \big\|^2 \\
   + (w-2)\int_{t_0}^{t} (s-1)^{w-3}
\big\| \A \tilde{X}(s) \big\|^2 ds
\bigg\} 
\le
\dfrac{w}{2}\left(\dfrac{w}{2}-1\right)
(t+1)^{w-2} \big\| \A \tilde{X}(t) \big\|^2 .
\end{multline}
We can drop the two positive terms on the left-hand side
(recall that $w \geq 2$), thus 
\begin{equation}\label{eq:Et-bd}
\E|_t
\le
\E|_{t_0}  +
\dfrac{w}{2}\left(\dfrac{w}{2}-1\right) (t+1)^{w - 2}
\big\| \A \tilde{X}(t)  \big\|^2 
\end{equation}
for all $t \ge t_0$.
Combining this with \eqref{eq:lyapunov4_2} (and ignoring the positive quadratic term)
we conclude that
\begin{equation}
\label{eq:boundphi1}
\alpha \tilde{\objective} \le \E|_{t_0} (t+1)^{-w} +
 w(w -2) \big(2 (t_0+1)\big)^{-2} \big\| \A \tilde{X}  \big\|^2
\end{equation}
for all $t \ge t_0$.
Since $\objective$ is strongly convex,
using a similar argument as that used to obtain \eqref{eq:stronga},  we have
\begin{equation}
\label{eq:strongaa}
\|\A \tilde{X} \|^2
\le  2 \sigma_1^2(\A)  \mu^{-1}  \tilde{\objective} .
\end{equation}
Using this inequality in the last term of \eqref{eq:boundphi1}, and
recalling the definition of $t_0$, yield
\begin{equation}
\label{eq:boundphi2}
\tilde{\objective} 
\le   2 \E|_{t_0} \alpha^{-1} (t+1)^{-w}   
\end{equation}
for all $t \ge t_0$.
Strong convexity of $\objective \equiv \objective_\lambda$
thus implies  
\begin{equation}
\big\| \tilde{X}_\lambda(t) \big\|^2  \le 2 \mu^{-1}  
\tilde{\objective} 
\leq  4 \E_\lambda|_{t_0}   (\alpha \mu)^{-1}    (t+1)^{-w }    
\end{equation}
for all $t \ge t_0$, where we restored the parameter $\lambda$.
Taking the limit $\lambda \downarrow 0$
yields \eqref{eq:relaxed_aadmm_strongly_convex} with
$C \equiv \E_{\lambda\downarrow 0}|_{t_0}$, i.e.,
\begin{equation}
\label{constant21}
C  =
(t_0 + 1)^w \big\{
 \alpha \left[ \objective(X(t_0)) - \objective^\star \right] 
+ (1/2) \big\| \A \big( w [ X(t_0) - x^\star ]  
+ (t_0+1) \dot{X}(t_0) \big) \big\|^2 \big\} ,
\end{equation}
with $t_0$ defined in \eqref{t0choice} and $w \equiv 2\rr / 3$.
\end{proof}

\begin{proof}[Proof of Proposition~\ref{th:heavy_ball_admm_convex}]
Consider the regularized differential equation \eqref{ode3}.
We again use the centered variables \eqref{eq:perturb}  and
omit the regularization parameter $\lambda$. Thus, consider
\begin{equation}
\label{eq:lyapunov5_2}
\E \equiv  \alpha t \tilde{\objective} 
+\dfrac{\rr}{2} \big\|  \A \tilde{X} \big\|^2
+ \dfrac{t}{2} \big\| \A \dot{\tilde{X}}\big\|^2 +
\big\langle \tilde{X} , \AtA \dot{\tilde{X}} \big\rangle .
\end{equation}
Its total time derivative along trajectories of \eqref{ode3} is
\begin{equation}
\begin{split}
\dot{\E} =
\alpha \big(
 \tilde{\objective}  - 
 \big\langle \tilde{X} , \nabla \tilde\objective
\big\rangle \big)
 -(t \rr - 1)\big\| \A \dot{\tilde{X}}\big\|^2 .
 \end{split}
\end{equation}
Convexity of $\objective$ tells us that the first term is
nonpositive, and so is the second term
provided $t \ge t_0 \equiv 1/\rr$. Therefore, under this condition, $\dot{\E} \le 0$.
We now proceed to show that $\E \ge 0$. Let us write \eqref{eq:lyapunov5_2} as
\begin{equation}
\label{eq:newlyap5}
\E =  \alpha t \tilde{\objective}  + \widehat{\E}
\end{equation}
with
\begin{equation}
\widehat{\E} \equiv \dfrac{\rr}{2}
\big\| \A \tilde{X} \big\|^2 +
\big\langle \tilde{X} , \AtA \dot{\tilde{X}} \big\rangle
+ \dfrac{t}{2}\big\|\A\dot{\tilde{X}}  \big\|^2.
\end{equation}
For all $t \geq t_0 = 1/\rr$ it follows that
\begin{equation}
\begin{split}
\widehat{\E} &\geq
\dfrac{\rr}{2} \big\| \A \tilde{X} \big\|^2 +
\big\langle \tilde{X} , \AtA \dot{\tilde{X}} \big\rangle  +
\dfrac{1}{2\rr}\big\|\A\dot{\tilde{X}} \big\|^2  \\
&= \dfrac{1}{2}\big\|\rr^{-1/2} 
\A\dot{\tilde{X}} +\rr^{1/2} \A \tilde{X} \big\|^2,
\end{split}
\end{equation}
showing that $\E \ge 0$.
We thus have
$\E(X(t), \dot{X}(t), t)  \le \E(X(t_0), \dot{X}(t_0), t_0)$ for all $t \ge t_0$, and
after dropping the positive term $\widehat{\E}$
we obtain
\begin{equation}
\objective_\lambda(X_\lambda(t)) - \objective^\star
\le  \dfrac{ \E_\lambda|_{t_0} }{ \alpha  t }  ,
\end{equation}
where we restored $\lambda$.
This implies \eqref{eq:heavy_ball_admm_convex} in the limit
$\lambda \downarrow 0$ with a constant 
\begin{equation}
\label{eq:const_hbadmm}
C \equiv  r^{-1} \alpha \left[ \objective(X(t_0)) - \objective^\star \right] \\
+ \dfrac{1}{2} \big\|  \rr^{1/2} \A \big( X(t_0) - x^\star \big)  + \rr^{-1/2} \A \dot{X}(t_0) \big\|^2. 
\end{equation}

For the strongly convex part, consider the function
\begin{equation}
\label{eq:lyapunov6_2}
\E \equiv  e^{2\rr t/3}
\Big\{   \alpha \tilde{\objective} 
+ (\rr/3)^2 \big\| \A \tilde{X}  \big\|^2
  + (1/2)  \big\| \A \dot{\tilde{X}} \big\|^2
+ ( 2\rr / 3)
\big\langle \A \tilde{X} , \A \dot{\tilde{X}} \big\rangle
\Big\} .
\end{equation}
Taking its total time derivative along trajectories of system
\eqref{ode3} yields 
\begin{equation}
\label{eq:lyapunov6_der}
\dot{\E} \le
\dfrac{2  \rr \alpha }{3} e^{2\rr t/3} 
\big( \tilde{\objective}  -
\big\langle \tilde{X}  , \nabla\tilde\objective \big\rangle
\big)  
+ \dfrac{2\rr^3}{27} e^{2\rr t /3} \big\|\A \tilde{X}  \big\|^2.
\end{equation}
Using the  inequality \eqref{eq:stronga} in the first term
yields
\begin{equation}
\dot{\E} \le \dfrac{\rr}{3} e^{2\rr t /3}
\left( \dfrac{2 \rr^2 }{9} -
\dfrac{\mu \alpha }{\sigma_1^2(\A)}\right)
\big\|\A \tilde{X} \big\|^2.
\end{equation}
Thus, if $\rr \le \bar{\rr}$, with $\bar{\rr}$ defined in
\eqref{eq:rbardef}, it follows that
$\dot{\E} \le 0$.
It remains to show that \eqref{eq:lyapunov6_2} is nonnegative. Note that this 
expression can be written as
\begin{equation} \label{new_E}
\E = (\alpha / 2) e^{2\rr t/3}
\tilde{\objective} 
+ e^{2\rr t/3}\widehat{\E}
\end{equation}
where
\begin{equation} \label{new_C_next}
\widehat{\E} \equiv
\dfrac{\alpha}{2} \tilde{\objective}
+\dfrac{\rr^2}{9} \big\|  \A \tilde{X}  \big\|^2   +
\dfrac{1}{2}\big\| \A \dot{\tilde{X}} \big\|^2 +
\dfrac{2\rr}{3}\big\langle \A \tilde{X}, \A\dot{\tilde{X}}
\big\rangle  .
\end{equation}
Using \eqref{eq:strongaa} and
defining
\begin{equation} \label{def:c12ab}
a \equiv \big\|\A \tilde{X} \big\|,  \qquad
b \equiv \big\| \A \dot{\tilde{X}} \big\|, \qquad
c_1 \equiv \dfrac{\mu \alpha}{4 \sigma_1^2(\A)} + \dfrac{r^2}{9},  \qquad
c_2 \equiv \dfrac{\rr}{3},  
\end{equation}
we have
\begin{equation}
\label{eq:EtildeSquare}
\begin{split}
\widehat{\E} &\ge
c_1 a^2  + \dfrac{1}{2} b^2   - 2 c_2 a b
 \\ &=   \left( \sqrt{c_1} a - \dfrac{c_2}{\sqrt{c_1}}
b \right)^2 +
\left( \dfrac{1}{2} -  \dfrac{c_2^2}{c_1} \right) b^2  .
\end{split}
\end{equation}
Thus, $\widehat{\E} \ge 0$,
since $c_2^2 \le c_1/2$ due to
$\rr \le \bar{\rr}$. This implies that $\E \ge 0$. Therefore, since $\dot{\E} \le 0$, we have
$\E(X(t), \dot{X}(t), t) \le \E(x_0, 0, 0)$, and upon neglecting
$\widehat{\E} \ge 0$ we get
\begin{equation}\label{diff-obj-last}
\objective(X(t)) - \objective^\star
\le   2 \E|_{t=0} \alpha^{-1}  e^{- 2 \rr t /3}
\end{equation}
provided $\rr \le \bar{\rr}$.
Also,
$\| \tilde{X} \|^2 \le (2/\mu)
\tilde{\objective} $ due to the
strong convexity of $\objective$, hence
\begin{equation} \label{Xtxstar}
\| X(t) - x^\star  \|^2 \le   4  \E\vert_{t=0} (\alpha \mu )^{-1} e^{-2 r t/3}
\end{equation}
for $\rr \le \bar{\rr}$.
Using \eqref{new_E}, \eqref{new_C_next},
\eqref{eq:strongaa},
and the condition $\rr \le \bar{r}$, we also have that
\begin{equation}
\begin{split}
\E\vert_{t=0}
&= \alpha \tilde{\objective}(x_0) 
+ (r^2 / 9)\| \A \tilde{x}_0  \|^2 \\
& \le \left( \alpha + 2 \rr^2 \sigma_1^2(A)/ (9 \mu)
   \right) \tilde{\objective}(x_0)  \\
& \le (3 \alpha / 2 ) \tilde{\objective}(x_0) .
\end{split}
\end{equation}
Finally, taking the limit $\lambda \downarrow 0$ 
into  \eqref{Xtxstar} (recall that we omitted the regularization parameter $\lambda$) yields
\eqref{eq:heavy_ball_admm_strongly_convex}
with a constant
$C = 6 \left( \objective(x_0) - \objective^\star\right)$.
\end{proof}

\subsection{Derivation and analysis of the modified equations}
\label{sec:deriv_mod}

In this part we assume  that $f$ and $g$ are smooth
(Assumption~\ref{th:assump4}) so the next derivation makes use of Taylor expansions
which are not available in the nonsmooth case.

\begin{proof}[Proof of Proposition~\ref{th:mod_eq_aadmm}]
From the optimality condition of updates 
\eqref{eq:raadmm2} and \eqref{eq:raadmm3} we conclude that
\begin{equation}\label{eq:ukg}
u_{k} = \epsilon^2 \nabla g(z_{k}) 
\end{equation}
for all $k=0,1,\dotsc$.
The optimality condition of \eqref{eq:raadmm1} reads
\begin{equation} \label{eq:optf}
\nabla f(x_{k+1}) + \A^T ( \A x_{k+1} - \hat{z}_k + \hat{u}_k ) / \epsilon^2 =0. 
\end{equation}
From the update \eqref{eq:raadmm4} and Eq.~\eqref{eq:ukg} we have
\begin{equation}\begin{split}
\hat{u}_k &= u_k + \gamma_k (u_{k} - u_{k-1}) \\ &= \epsilon^2 \nabla g(z_k) + \gamma_k
\epsilon^2\big(\nabla g(z_k) - \nabla g(z_{k-1})\big) .
\end{split}
\end{equation}
Note that for the R-A-ADMM method \eqref{eq:relaxed_aadmm} we have
\begin{equation} \begin{split} \label{eq:gammat1}
\gamma(t) 
= 1 - \dfrac{\rr \epsilon}{t+1} + \dfrac{\rr (\rr - 1) \epsilon^2}{(t+1)^2} +\bigO(\epsilon^3) ,
\end{split}
\end{equation}
so that  $\hat{u}_k \to \hat{U}(t)$ as $\epsilon \to 0$, where
\begin{equation} \label{eq:uhat4}
\hat{U}(t) = \epsilon^2 \nabla g(Z(t)) + \epsilon^3 \nabla^2 g(Z(t)) \dot{Z}(t) +\bigO(\epsilon^4).
\end{equation}
On the other hand, for the R-HB-ADMM method \eqref{eq:heavy_ball_admm} we have
\begin{equation}
\gamma(t) = 1 - \rr \epsilon ,
\end{equation}
which yields the same conclusion \eqref{eq:uhat4}.  Thus, from Eq.~\eqref{eq:optf}, 
for both methods \eqref{eq:relaxed_aadmm} and \eqref{eq:heavy_ball_admm} it holds that
\begin{equation} \label{eq:optf2}
\nabla f(X(t+\epsilon)) + 
\A^T \nabla g(Z(t)) + \epsilon \A^T \nabla^2 g(Z(t)) \dot{Z}(t)  + 
\A^T \dfrac{ \A X(t+\epsilon) - \hat{Z}(t) }{\epsilon^2} =
 \bigO(\epsilon^2) .
\end{equation}

Consider the update \eqref{eq:raadmm3}, i.e.,
\begin{equation} \label{eq:axhz}
\A x_{k+1} - \hat{z}_k = \alpha^{-1}\left( u_{k+1} - \hat{u}_k + z_{k+1} - \hat{z}_k \right).
\end{equation}
In light of \eqref{eq:ukg} and \eqref{eq:uhat4} we conclude that
\begin{equation}\label{eq:uuzero}
u_{k+1} - \hat{u}_k \to U(t+\epsilon) - \hat{U}(t) = \bigO(\epsilon^4) ,
\end{equation} 
so this term can be neglected; it contributes to $\bigO(\epsilon^2)$ in \eqref{eq:optf2}.  Using update \eqref{eq:raadmm5} and \eqref{eq:gammat1},
 $z_{k+1}-\hat{z}_k$ tends to
\begin{equation}\label{eq:zhz1}
Z(t+\epsilon) - \hat{Z}(t) 
 = \epsilon^2 \left[ \ddot{Z}(t) + \dfrac{\rr}{t+1} \dot{Z}(t) \right] 
 - \epsilon^3 \left[ \dfrac{r}{2(t+1)} \ddot{Z}(t) + \dfrac{r(r-1)}{(1+t)^2} \dot{Z}(t)  \right] + \bigO(\epsilon^4).
\end{equation}
On the other hand, with \eqref{eq:heavy_ball_admm5} we conclude that
\begin{equation}\label{eq:zhz2}
z_{k+1} - \hat{z}_k \to \epsilon^2\left[ \ddot{Z}(t) + \rr \dot{Z}(t)  \right]
- (\epsilon^3 \rr / 2) \ddot{Z}(t) +  \bigO(\epsilon^4).
\end{equation}
Collecting these last steps, we have from \eqref{eq:optf2} that
the continuous-time limit of R-A-ADMM 
\eqref{eq:relaxed_aadmm}  is given by
\begin{equation}\label{eq:eq_cont_aadmm} \begin{split}
& \nabla f(X(t)) + \A^T \nabla g(Z(t)) 
+ \epsilon \nabla^2 f(X(t)) \dot{X}(t) + \epsilon \A^T \nabla^2 g(Z(t)) \dot{Z}(t) \\
&+ \alpha^{-1}\A^T\bigg[  \ddot{Z}(t) + \dfrac{\rr}{t+1} \dot{Z}(t) \bigg]  
- \alpha^{-1}\A^T\bigg[ \dfrac{\epsilon r}{2(t+1)} \ddot{Z}(t) + \dfrac{\epsilon r(r-1)}{(1+t)^2} \dot{Z}(t) 
\bigg]  = \bigO(\epsilon^2) ,
\end{split}
\end{equation}
while for R-HB-ADMM \eqref{eq:heavy_ball_admm} we have
\begin{equation} \label{eq:eq_cont_hb_admm} \begin{split}
& \nabla f(X(t)) + \A^T \nabla g(Z(t)) 
+ \epsilon \nabla^2 f(X(t)) \dot{X}(t) + \epsilon \A^T \nabla^2 g(Z(t)) \dot{Z}(t) \\
&+ \alpha^{-1}\A^T\big[  \ddot{Z}(t) +  \rr \dot{Z}(t) - (\epsilon r / 2 ) \ddot{Z}(t)\big]  
  = \bigO(\epsilon^2) .
\end{split}
\end{equation}

We need a relation between $X(t)$ and $Z(t)$.
From \eqref{eq:axhz}--
\eqref{eq:zhz2} we have
$\A X(t+\epsilon) - \hat{Z}(t) = \bigO(\epsilon^2)$.
Using \eqref{eq:raadmm5} and \eqref{eq:gammat1} this last equation yields
$Z(t) = \A X(t)  + \epsilon \big(\A \dot{X}(t) - \dot{Z}(t)\big) = \bigO(\epsilon^2)$,
thus differentiating this expression and replacing the result into the second term yield
\begin{equation}
Z(t) = \A X(t) + \bigO(\epsilon^2).
\end{equation}
Using this into \eqref{eq:eq_cont_aadmm} we obtain  the continuous-time limit of algorithm \eqref{eq:relaxed_aadmm}, up to an approximation of 
$\bigO(\epsilon^2)$, 
as
\begin{equation} \label{eq:mod_eq_aadmm2}
\alpha^{-1}\AtA \left[ \beta(t) \ddot{X} + \gamma(t) \dot{X}  \right]
+ \epsilon \nabla^2 \objective(X) \dot{X} = - \nabla \objective(X) 
\end{equation}
with
\begin{equation} \label{eq:coef_mod_aadmm}
\beta(t) \equiv 1 - \dfrac{\epsilon \rr}{2(t+1)}, \qquad 
\gamma(t) \equiv \dfrac{r}{t+1} - \dfrac{\epsilon r (\rr-1)}{(t+1)^2}.
\end{equation}
Note that used $\nabla \objective(X) = \nabla f(X) + \A^T \nabla g(\A X)$ and 
$\nabla^2 \objective(X) = \nabla^2 f(X) + \A^T \nabla^2 g(\A X) \A$ in obtaining
this equation.
Similarly, using \eqref{eq:eq_cont_hb_admm}
we obtain for algorithm \eqref{eq:heavy_ball_admm} the same differential equation
\eqref{eq:mod_eq_aadmm2} but with coefficients
\begin{equation} \label{eq:coef_hb_admm}
\beta \equiv 1 - \dfrac{ \epsilon \rr}{2}, \qquad \gamma(t) \equiv \rr .
\end{equation}
Now, consider dividing equation
\eqref{eq:mod_eq_aadmm2} by $\beta(t)$ and expanding the coefficients up to 
$\bigO(\epsilon^2)$; this recasts \eqref{eq:mod_eq_aadmm2} into
\eqref{eq:mod_eq_aadmm3}. 
\end{proof}

Next, we show that all previous rates for the leading order dynamics 
are preserved for the perturbed system \eqref{eq:mod_eq_aadmm3}.

\begin{proof}[Proof of Proposition~\ref{th:preserve_rates_pert}]
Recall the  function~\eqref{eq:lyap_pert}. 
We showed in \eqref{eq:der_lyap_pert} that this function 
is decreasing, i.e.,\footnote{In the case of decaying damping \eqref{eq:ab_dec} it
is necessary that $a(t) > 0$, which implies  $\epsilon < 2/(r-2)$.
For $\rr  = 3$ this requires the step size to obey $\epsilon < 2$. We assume a suitable
step size from now on. Similar requirement is not necessary for the constant
damping case \eqref{eq:ab_const}.}
\begin{equation}\label{eq:dotE_pert}
\dot{\E} \le -a(t) \alpha^{-1} \big\| \A \dot{X} \big\|^2 
- \epsilon \big\langle \dot{X}, \nabla^2 \objective \dot{X} \big\rangle \le 0.
\end{equation}
When $\epsilon \to 0$ this recovers the similar result for the leading
order systems \eqref{eq:relaxed_aadmm_ode} and \eqref{eq:relaxed_heavy_ball_aadmm_ode},
 which has less contraction due to the absence of the
Hessian term above. 
 Denote by $X_\epsilon \equiv X_\epsilon(t)$ the trajectory of the perturbed 
 system~\eqref{eq:mod_eq_aadmm3}, and by $X \equiv X(t)$
the  trajectory of the corresponding leading order system (with $\epsilon \to 0$).
The inequality \eqref{eq:dotE_pert} shows that
\begin{equation}
\dot{\E}(X_\epsilon, \dot{X}_\epsilon, t) \le \dot{\E}(X, \dot{X}, t),
\end{equation}
which upon integrating both sides yields
\begin{equation} \label{eq:bound_pert}
\E(X_\epsilon, \dot{X}_\epsilon, t) \le \E(X, \dot{X}, t),
\end{equation}
since $X_\epsilon(0) = X(0)$ and $\dot{X}_\epsilon(0) = \dot{X}(0) = 0$ (both systems
have the same starting point).
Thus, trajectories of the perturbed system are always bounded by  trajectories
of the leading order system.  Now, we already know convergence rates
for $\objective(X(t)) - \objective^\star$.
If we obtain a similar rate for $\E(X, \dot{X}, t)$, i.e., for the term $\| \A \dot{X}\|^2$,
then we can immediately conclude that the same rate holds 
for $\objective(X_\epsilon(t)) - \objective^\star$.
Next, we show that $\big\| \A \dot{X} \|^2$ decays at the same rate
as $\objective(X(t)) -\objective^\star$.

Write both systems \eqref{eq:relaxed_aadmm_ode} and \eqref{eq:relaxed_heavy_ball_aadmm_ode}
(in the smooth setting) as 
\begin{equation} \label{eq:eta_eq_form}
\alpha^{-1}\AtA\left[ \ddot{X} + \dot{\eta}(t) \dot{X}  \right] = -\nabla \objective(X),
\end{equation}
where $\eta(t) = \rr \log(t+1)$ for system \eqref{eq:relaxed_aadmm_ode} and
$\eta(t) = \rr t$ for system \eqref{eq:relaxed_heavy_ball_aadmm_ode}.
Letting
\begin{equation}
h(t) \equiv \big\| \A \dot{X}(t)\big\|^2
\end{equation}
and taking the inner product of \eqref{eq:eta_eq_form} with $\dot{X}$ we obtain
\begin{equation}
\dfrac{d}{dt}\left( e^{2\eta(t)} h(t)  \right) = -2 \alpha e^{2\eta(t)} \dfrac{d}{dt}
\left( \objective(X(t)) - \objective^\star\right).
\end{equation}
Thus, upon integrating both sides we get
\begin{equation}
 h(t) = -2 \alpha e^{-2\eta(t)}\int_0^t ds \, e^{2\eta(s)} \dfrac{d}{ds}\left( \objective(X(s)) - \objective^\star\right) .
\end{equation}
Integration by parts yields
\begin{equation} \label{eq:sol_h}
\begin{split}
h(t) &= 
%
-2 \alpha \left( \objective(X(t)) - \objective^\star\right)
+ 2\alpha e^{-2\eta(t)} \left( \objective(X(0)) - \objective^\star \right) \\
&\qquad + 2 \alpha e^{-2\eta(t)} \int_0^t ds \, e^{2\eta(s)}\dot{\eta}(s) \left( \objective(X(s)) - \objective^\star \right) .
\end{split}
\end{equation}

Consider the rate \eqref{eq:relaxed_aadmm_convex} where $\eta(t) = \rr \log(t+1)$.
Then the term with the integral in \eqref{eq:sol_h} is less or equal than
\begin{equation}
2\alpha (t+1)^{-2\rr} C \rr \int_0^t ds \, (s+1)^{2\rr - 3} =
\bigO\left( (t+1)^{-2} \right) .
\end{equation}
The 1st term in \eqref{eq:sol_h} also decays as $\sim (1+t)^{-2}$, and
the 2nd term decays faster, as $\sim (t+1)^{-2\rr} $.
Therefore,
\begin{equation}
h(t) = \big\| \A \dot{X}(t)\big\| = \bigO\left( (1+t)^{-2 }\right) .
\end{equation}
Consequently, $\E(X, \dot{X}, t) = \bigO\big((t+1)^{-2}\big)$ for the
leading order system. From \eqref{eq:bound_pert} the same holds true for the
perturbed system. From the definition of $\E$ in \eqref{eq:lyap_pert} we have
$\objective(X_\epsilon(t)) - \objective^\star \le \E$, thus
there exists some constant $K > 0$ such that 
\begin{equation}
\objective\big(X_\epsilon(t) \big) - \objective^\star \le \dfrac{K }{(t+1)^{2} } ,
\end{equation}
i.e., the perturbed system has at least the same convergence rate as the original
leading order system.

Consider now the convergence rate \eqref{eq:heavy_ball_admm_strongly_convex}
where $\eta(t) =  \rr t$ for the system with constant damping.
In terms of $\objective$ (see \eqref{diff-obj-last}) we have 
$\objective(X(t)) - \objective^\star \le C e^{-2\rr t/3}$.
The term with the integral in \eqref{eq:sol_h} becomes
\begin{equation}
2\alpha e^{-2 \rr t} C \rr \int_0^t ds \, e^{4\rr s/3} = \bigO\big( e^{-2\rr t/3}\big) .
\end{equation}
The 1st term with the objective function in \eqref{eq:sol_h} has the same decay rate,
while the 2nd term decays faster as $\sim e^{-2 \rr t}$. Therefore 
$h(t) = \bigO\big(e^{-2\rr t/3}\big)$ and also
$\E(X_\epsilon, \dot{X}_\epsilon) = \bigO\big(e^{-2\rr t/3}\big)$. 
From \eqref{eq:bound_pert} we thus conclude that
\begin{equation}
\objective\big(X_\epsilon(t)\big) - \objective^\star = \bigO\big(e^{-2\rr t/3} \big).
\end{equation}
By the same argument from \eqref{diff-obj-last} to \eqref{Xtxstar}, i.e., using
strong convexity of $\objective$, we also have  
$\| X_\epsilon(t) - x^\star \|^2 = \bigO\big(e^{-2\rr t/3}\big)$ so that
\eqref{eq:heavy_ball_admm_strongly_convex} holds true
for the perturbed system \eqref{eq:mod_eq_aadmm3}/\eqref{eq:ab_const}.

The exact same argument applies  
to the convergence rates \eqref{eq:relaxed_aadmm_strongly_convex}
and \eqref{eq:heavy_ball_admm_convex} for the respective perturbed systems as well.
\end{proof}

\section{Numerical experiments}
\label{sec:numerical}

\subsection{Trends in time series}
Consider the problem of estimating piecewise linear trends in time series data.
This can be done by solving \cite{Boyd:2009}
\begin{equation}
\label{eq:l1trend}
\min_{x,z} \dfrac{1}{2}\| y - x\|^2  +
\lambda \| z \|_1 \quad \mbox{such that} \quad z = \bm{D} x ,
\end{equation}
where $y \in \mathbb{R}^n$ is a given signal,
$\bm{D} \in \mathbb{R}^{(n-2)\times n}$ is a Toeplitz
matrix with first row $(1,-2,1,0,\dotsc,0)$, and $\lambda > 0$ is the
regularization parameter. (For this type of problem $\lambda$ has to
be large, e.g., $\lambda \ge \| (\bm{D}\bm{D}^T)^{-1} \bm{D} y\|_{\infty}$; see \cite{Boyd:2009}.)
The above problem is well-suited to the variants of ADMM since 
the proximal operators of $f(x)$ and $g(z)$
have a closed form solution \cite{Boyd:2011}.
Consider a time series
\begin{equation} \label{eq:time_series}
y_i = x_i + \xi_i, \qquad x_{i+1} = x_i + v_i, \qquad x_0 = 0,
\end{equation}
for $i=1,\dotsc,n$. Here $x\equiv (x_1, \dotsc, x_n)$ is the true underlying
trend, which is superimposed by noise $\xi_i \sim \mathcal{N}(0, \sigma^2)$.
The  slopes are generated by a Markov process where $v_{i+1} = v_i$ with
probability $p$ and
$v_i \sim \mathcal{U}(-b,b)$ with probability $(1-p)$ for some $b > 0$.
The goal is to recover $x$ from $y$.

\begin{figure}[t]
\centering
\includegraphics[scale=.43]{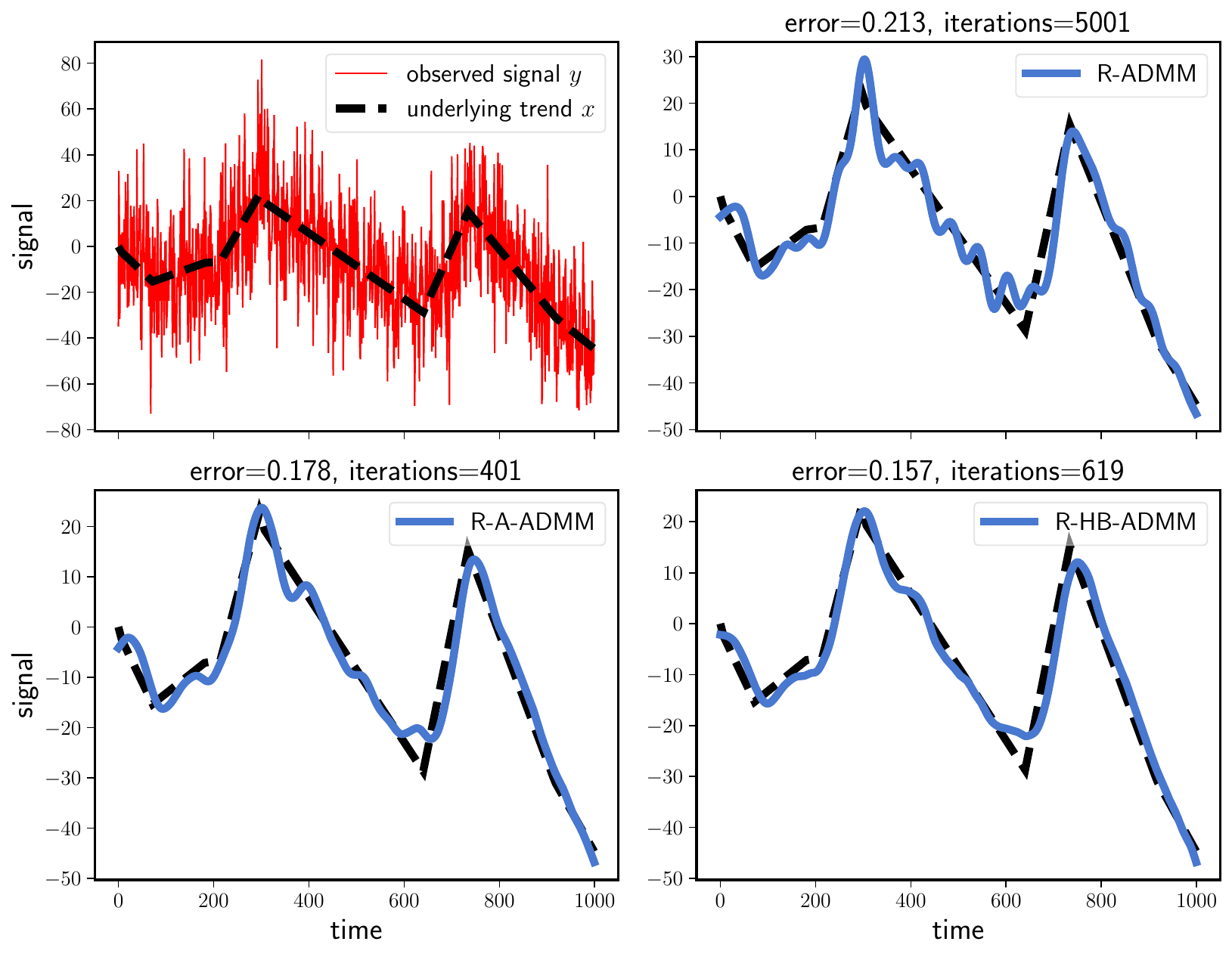}
\caption{Reconstructing the
piecewise linear trend $x$ in \eqref{eq:time_series} (black dashed line) by
solving  \eqref{eq:l1trend}.
We indicate the relative error,
$\| \hat{x} - x \|/\| x \|$, where $\hat{x}$ is the algorithm's
estimate of $x$, as well as the number of iterations.}
\label{fig:recovered_trend}
\end{figure}

\begin{figure}
\centering
\includegraphics[scale=0.5]{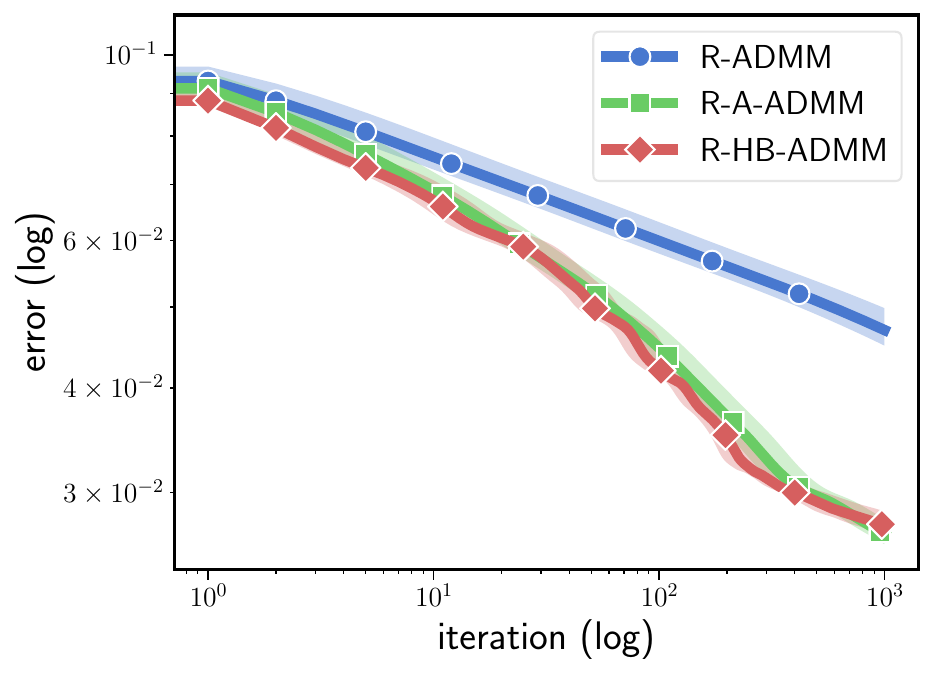}\qquad
\includegraphics[scale=0.5]{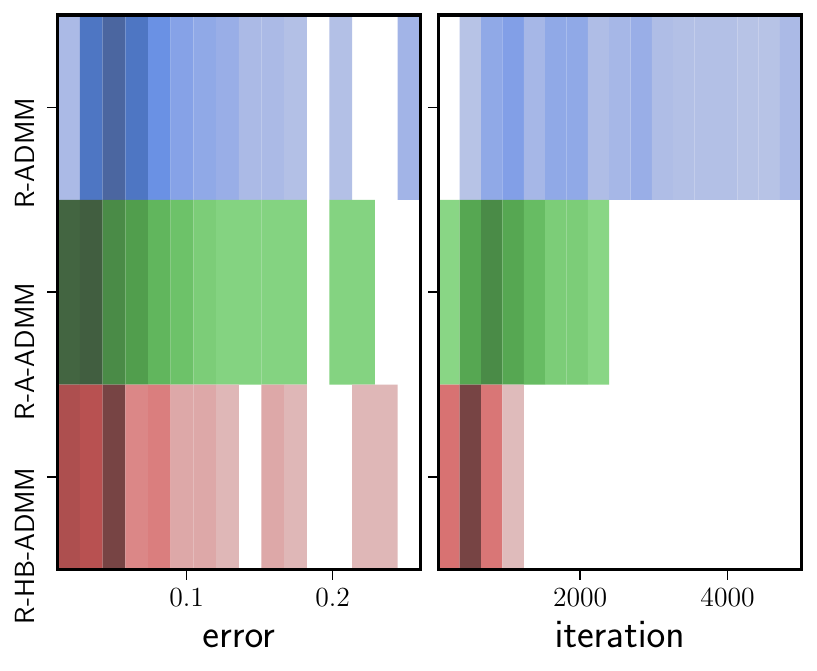}
\caption{\emph{Left:} Convergence rate for an instance of the problem
in Fig.~\ref{fig:recovered_trend}.
Solid lines correspond to $\alpha=1$ (relaxation parameter of ADMM) and the
shaded areas to $\alpha \in [0.6, 1.4]$ ($\alpha > 1$ is faster).
\emph{Right:} Histograms with relative error and number
of iterations over 80 random trials. The algorithms run until a small
tolerance condition is achieved.}
\label{fig:trend_plots}
\end{figure}

We consider the above model with
$n=1000$, $p=0.99$, $\sigma = 20$, and $b=0.5$.
We choose $\lambda = 15 \times 10^3$ in \eqref{eq:l1trend} and
inverse step size $\rho = 1$ for all algorithms.
For R-A-ADMM~\eqref{eq:relaxed_aadmm} we use
the default $\rr = 3$.
For R-HB-ADMM~\eqref{eq:heavy_ball_admm} we choose  $\gamma = 0.01$.
Fig.~\ref{fig:recovered_trend} shows the recovered trend
for one sample of this model. Both accelerated variants are
faster and more accurate than the base method.
Fig.~\ref{fig:trend_plots} shows the convergence rates
for another sample of this problem, and 
histograms of the relative error and number of iterations
over 80 Monte Carlo runs.
 Both accelerated variants improve over
R-ADMM, and  R-HB-ADMM converged in less iterations
than R-A-ADMM.

\subsection{Robust principal component analysis}
Let
$\M = \X^\star + \Z^\star \in \mathbb{R}^{n\times m}$,
where $\X^\star$ has low rank and
$\Z^\star$ is sparse. Under certain rank and sparsity conditions
it is possible to recover 
$\X^\star$ and $\Z^\star$ from
observation of $\M$ alone. This is done by solving \cite{Candes:2011}
\begin{equation} \label{eq:pcp}
  \min_{\X,\Z} \| \X\|_* + \lambda \| \Z \|_1 \quad \mbox{such that}
  \quad \X + \Z = \M ,
\end{equation}
where $\lambda = 1/\max\{n ,m\}$ and 
$\| \X \|_*$ 
is the nuclear norm. 
This problem
is known as \emph{robust principal component analysis} (robust PCA)
and can be seen as an idealized version of PCA for highly corrupted data.
PCA is arguably one of the most used techniques for dimensionality reduction
\cite{BishopBook,Jolliffe:2016}.
(Robust) PCA has
several important applications in statistics, signal processing,
and machine learning \cite{Zou:2006,Jordan:PCA,Pados:2017}.

\begin{figure}[t]
\centering
\includegraphics[scale=.45]{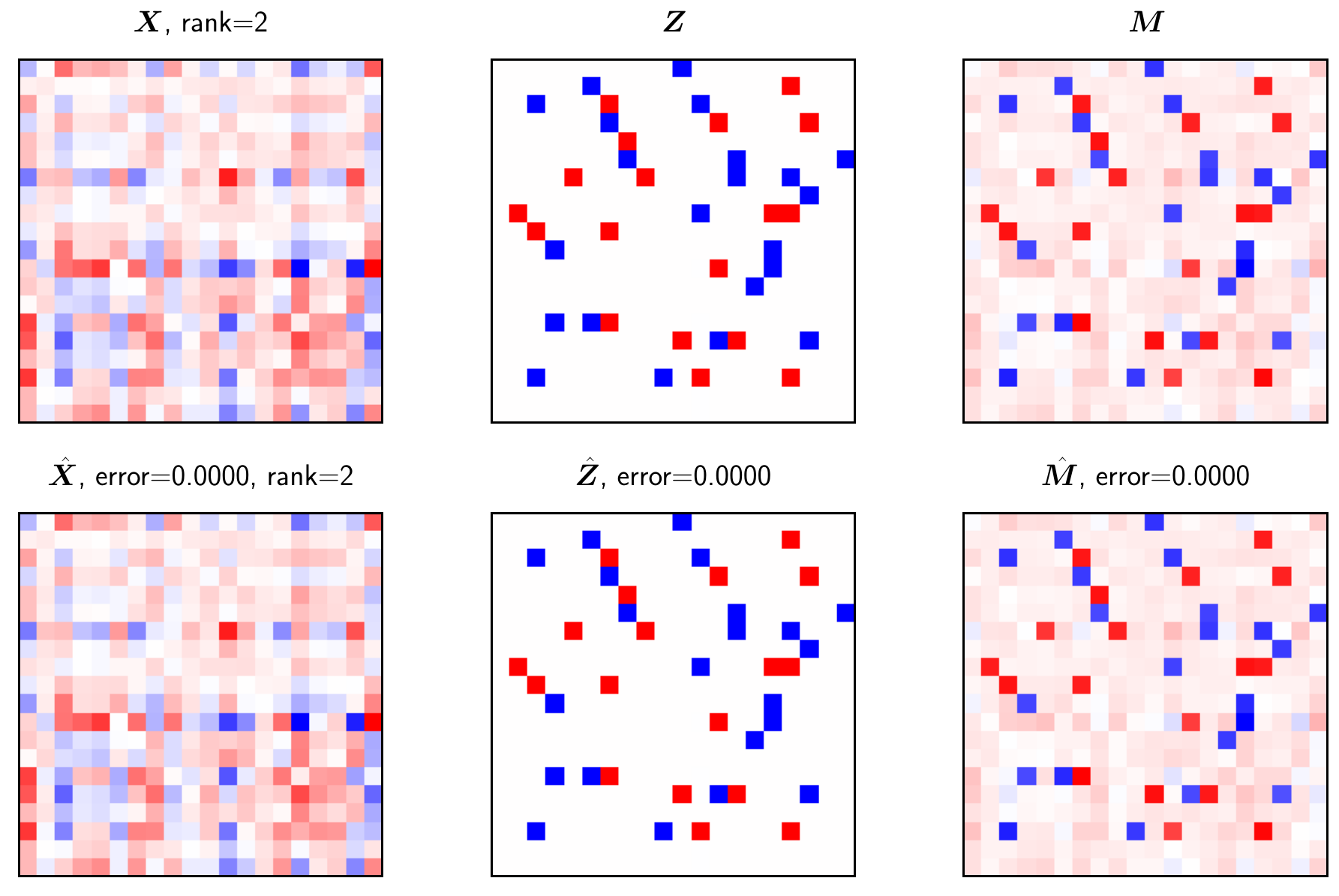}
\put(-265,180){$\bm{+}$}
\put(-130,180){$\bm{=}$}
\put(-265,60){$\bm{+}$}
\put(-130,60){$\bm{=}$}
\caption{Recovered matrices when solving \eqref{eq:pcp}.
We show the results for H-HB-ADMM. The relative
error is
$\| \hat{\bm{X}} - \bm{X}\|/\| \bm{X}\|$ where $\hat{\bm{X}}$
is the algorithm's estimate of $\bm{X}$, and analogously for
the other matrices. Note that both $\X,\Z$ are recovered exactly.}
\label{fig:mat}
\end{figure}

\begin{figure}[t]
\centering
\includegraphics[scale=.5]{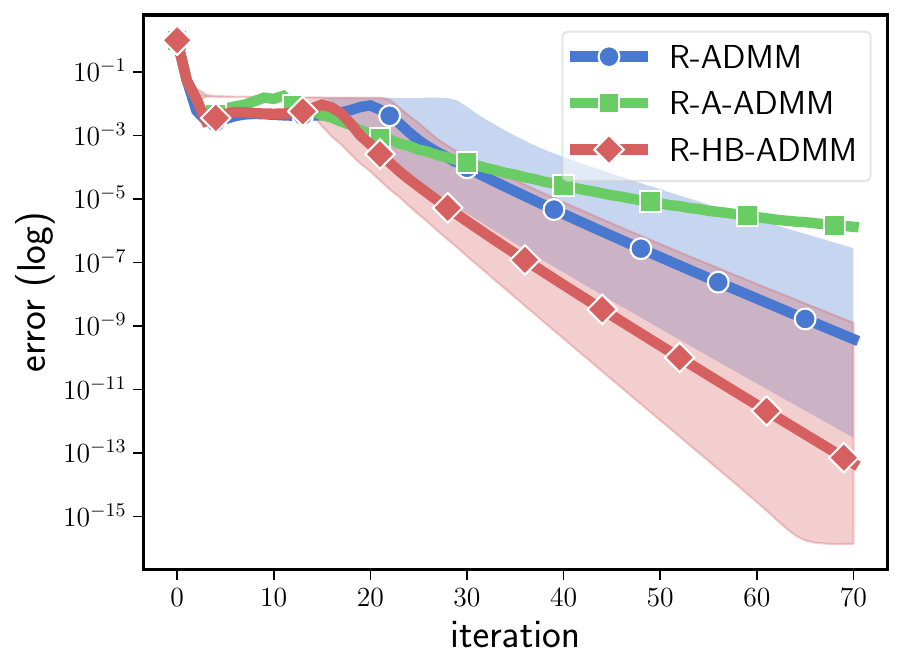}\qquad
\includegraphics[scale=.5]{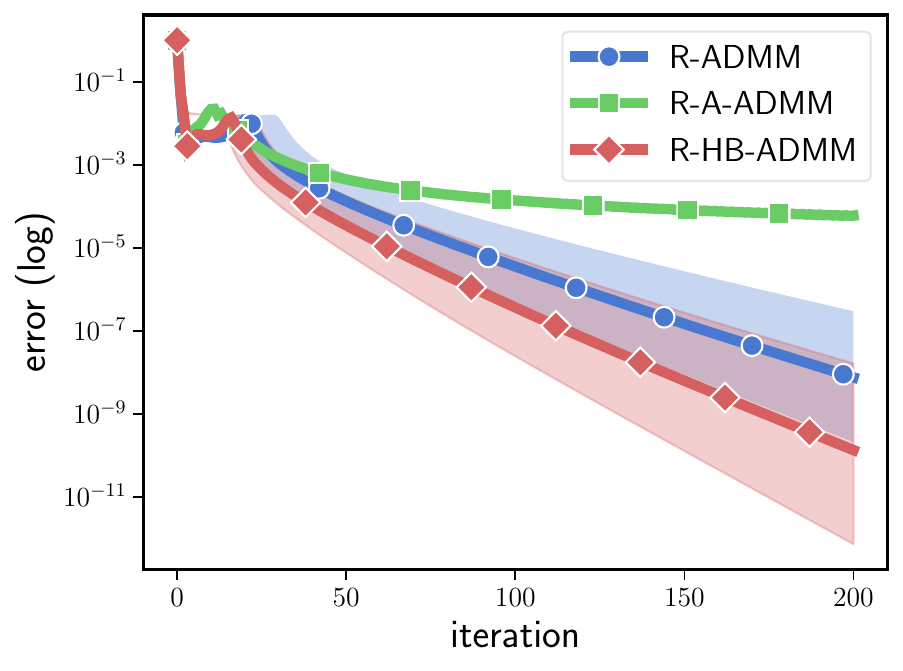}
\caption{
Robust PCA \eqref{eq:pcp} where exact recovery is possible. 
We report the relative error
$\| \X_k + \Z_k - \M\|/\| \M\|$ versus the iteration $k$.
Solid lines correspond to $\alpha=1$ and the shaded
areas $\alpha\in [0.7, 1.3]$.
\emph{Left:} in this case exact recovery is guaranteed.
\emph{Right:} this is a more challenging case, on the boundary of
the phase transition where exact recovery is impossible.
($\alpha\ne 1$ for R-A-ADMM did not improve over $\alpha=1$ and this is why
there is no shaded area in this case.)}
\label{fig:rpca}
\end{figure}

It has been noted \cite{Candes:2011} that ADMM
with  $\rho=1$
is very effective in solving \eqref{eq:pcp}, faster and more
accurate than several methods including nonsmooth extensions of Nesterov's
method. Thus,
we wish to verify whether our accelerated variants of ADMM
are able to improve over standard ADMM.
The proximal operators for $\| \X \|_*$ and $\| \Z\|_1$
have well-known closed form expressions \cite{Candes:2011}.
We generate a matrix $\M \in \mathbb{R}^{n\times n}$ with
$\X^\star = \M_1 \M_2^T$, where
$\M_{1,2}  \sim \mathcal{N}(0, 1/n)$ are $n\times q$ matrices and
$\Z^\star \in \{ -1, 0, 1 \}^{n\times n}$ has Bernoulli $\pm 1$ entries, and
support of size $s$ chosen uniformly at random.
The results in a setting where exact
recovery is possible can be obtained by setting
$q = 0.05 \times n$ and $s = 0.1 \times n^2$
(this was considered in \cite{Candes:2011}). For all algorithms,
we use $\rho = 1$. For R-A-ADMM we use $r=3$, and for H-HB-ADMM we choose
$\gamma = 0.75$.  For illustration purposes, in Fig.~\ref{fig:mat} we show
a simple case with $n=20$ where R-HB-ADMM recovers each component
exactly.
In Fig.~\ref{fig:rpca} (left) we use the same setting but
in higher dimensions,  $n=1000$, and show the convergence rate of each algorithm.
A more challenging case,  close to the phase transition boundary
where exact recovery fails, is shown in Fig.~\ref{fig:rpca} (right);
the only difference is that we set a higher rank $q = 0.2 \times n$.
Interestingly, R-HB-ADMM
improves over R-ADMM, whereas R-A-ADMM does not.
The improvement of R-ADMM and R-HB-ADMM with $\alpha > 1$ agrees
with our theoretical results.
In these examples, choosing $\alpha > 1$ for R-A-ADMM did not improve
over $\alpha=1$.

Although problem \eqref{eq:pcp} is convex, 
the plots in Fig.~\ref{fig:rpca}
show that R-ADMM and R-HB-ADMM have approximately linear convergence for
sufficiently large $k$.
This suggests that
there is a region in which the objective function behaves as a strongly
convex function. R-HB-ADMM has  improved convergence
over R-ADMM,
as predicted by Proposition~\ref{th:heavy_ball_admm_convex}.
On the other hand,
Proposition~\ref{th:relaxed_aadmm_convex} tells us that R-A-ADMM is
unlikely to attain linear
convergence, as  reflected in Fig.~\ref{fig:rpca}.
These empirical results are 
consistent with our theoretical findings.

\section{Conclusion}
\label{sec:discussion}

We introduced two new families of \emph{relaxed and
accelerated} ADMM algorithms
and derived differential inclusions
modeling these
methods in the limit of small step size. 
We  obtained convergence rates for these systems in the convex and
strongly convex cases, highlighting an interesting tradeoff
between the type of dissipation  versus degree of convexity of the objective function; the
complexity results are summarized
in Table~\ref{table:convergence}.
These rates remain unproven for the algorithms,
whose discrete analyses are challenging.   However, we provided backward error
analysis by deriving modified or perturbed
differential 
equations that model these algorithms more closely, even for finite step sizes, for which
such convergence rates hold true and the perturbed system has improved stability.

Due to the nonsmooth nature of the differential inclusions, such systems  are
able to model the variants of ADMM when applied to some types of
constrained  problems. In Appendix~\ref{sec:constrained} 
we provide
further details on how to generalize these methods to account 
for more general  nonlinear constraints, and make connections with
multiple scale analysis and singular perturbation theory. 

Let us make one final remark about 
the proofs of convergence rates. 
A closer inspection shows that only a weaker notion
of (strong) convexity  was 
actually needed, i.e., it is sufficient
to consider the relation only at the critical point,
\begin{equation}\label{eq:quasi_strong}
\objective(x^\star) \ge \objective(x)
+ \langle \partial \objective(x), x^\star -  x \rangle
+ (\mu/2) \| x - x^\star\|^2
\end{equation}
for all $x \in \dom\objective$.
This can be taken as the definition
of \emph{quasi-strong convexity} (\emph{quasi-convexity} when $\mu=0$) \cite{Necoara:2019}.
Thus, we expect that the previous convergence rates  hold  beyond
(strong) convexity conditions, and perhaps even for some  
nonconvex problems.\footnote{%
In particular, we expect that the exponential
convergence \eqref{eq:heavy_ball_admm_strongly_convex}  remains valid
under  \eqref{eq:quasi_strong}, which would explain the success of the method with
constant damping for nonstrongly convex problems, such as the one in Fig.~\ref{fig:rpca}.}
When  $\objective$ is smooth and $\nabla \objective$ is
Lipschitz continuous\footnote{Which is a sufficient but not 
necessary for uniqueness of solutions.}   
this  follows straightforwardly.
However, when $\objective$ is nonsmooth the situation is 
more delicate
because \eqref{eq:quasi_strong} does not even imply convexity, thus 
uniqueness of resolvents, Moreau-Yosida regularization,
and solutions of differential inclusions become subtle issues. 


\subsection*{Acknowledgements}
This work was supported by grants ARO MURI W911NF-17-1-0304, 
NSF 2031985, and NSF 1934931. 
We thank the anonymous referees for insightful comments.

\appendix

\section{Derivation of the differential inclusions}
\label{sec:der_inclusions}

In this section we provide the  derivations of the leading order 
differential inclusions modeling the 
ADMM algorithms introduced in Section~\ref{sec:variants}.

\begin{proof}[Proof of Proposition~\ref{th:relaxed_admm_ode}]
The optimality conditions for the updates \eqref{eq:relaxed_admm1} and \eqref{eq:relaxed_admm2}
read
\begin{align}
0 &\in \epsilon \partial f(x_{k+1}) + \A^T \left( \A x_{k+1} - z_k + u_k \right) , \\
0 &\in \epsilon  \partial g(z_{k+1}) - \alpha \A x_{k+1} - (1-\alpha) z_k + z_{k+1} - u_k .
\end{align}
Making use of \eqref{eq:relaxed_admm3} into the second equation above gives
\begin{equation} \label{eq:uk}
u_{k} \in \epsilon \partial g(z_{k})  ,
\end{equation}
which replaced into the first equation yields
\begin{equation}
\A^T \left( \A x_{k+1} - z_k  \right) \in - \epsilon \partial f(x_{k+1}) - \epsilon \A^T\partial g(z_k) .
\end{equation}
Using again the update \eqref{eq:relaxed_admm3} we can write this as
\begin{equation} \label{eq:comb}
\alpha^{-1} \A^T \left(  u_{k+1} - u_k + z_{k+1} - z_k \right)  \in -\epsilon \partial f(x_{k+1}) - \epsilon \A^T \partial g(z_k).
\end{equation}
Recall the discussion following Assumption~\ref{assump2}.
We have a discrete grid $t_k = \epsilon k$ and
there exists continuous functions $X, Z, U$ such that
$(X(t_k), Z(t_k), U(t_k)) \approx (x_k, z_k, u_k)$, where this approximation
is exact when $\epsilon \to 0$.
Henceforth we omit the $k$ dependence, i.e., $t \equiv t_k$,
anticipating the $\epsilon \to 0$ limit.
From Eq.~\eqref{eq:uk} we  have
$U(t) \in \epsilon \partial g(Z(t)) \to 0$ as $\epsilon \to 0$, and 
also $\dot{U}(t) = 0$ as $\epsilon \to 0$. 
In the limit $\epsilon\to 0$, Eq.~\eqref{eq:relaxed_admm3} implies
$Z(t) = \A X(t)$, and thus also $\dot{Z}(t) = \A \dot{X}(t)$.
Relation  \eqref{eq:comb} thus becomes
\begin{equation}
\alpha^{-1} \A^T \left( \dfrac{U(t+\epsilon)-U(t)}{\epsilon} + \A \dfrac{X(t+\epsilon) - X(t)}{\epsilon}  \right) 
\in -\partial f(X(t+\epsilon)) - \A^T \partial g(\A X(t)) .
\end{equation}
Taking the limit $\epsilon \to 0$,
and recalling that $\dot{U}(t) = 0$,  
yield the differential inclusion
\begin{equation}
\alpha^{-1} \AtA \dot{X}(t) \in - \partial f(X(t)) - \A^T \partial g(\A X(t)). 
\end{equation}
From \eqref{lemma}
this yields 
\eqref{eq:relaxed_admm_ode}.
This is a first-order system so the
dynamics is specified by the initial condition $X(0)=x_0$.
\end{proof}

\begin{remark} \label{rmk:initial}
We mention a subtlety regarding the initial condition in
Proposition~\ref{th:relaxed_admm_ode} which  also applies to the other
differential inclusions obtained in this paper.
It is necessary that $X(0) = x_0$ matches the starting point of the
algorithm.
Recall that \eqref{eq:uk} implies $U(t) = 0$
($\epsilon \to 0$) for any $t \in \mathcal{I}$.
We thus assume that \eqref{eq:relaxed_admm} is initialized
with some $x_0 \in \dom \objective$, $z_0 = \A x_0$, and $u_{0} = 0$;
this initialization is also assumed in algorithms \eqref{eq:relaxed_aadmm}
and \eqref{eq:heavy_ball_admm} (and moreover $\hat{u}_{0} = u_0 =  0$
and $\hat{z}_{0} = z_0 = \A x_0$ for the so-called  ``accelerated variables'').
\end{remark}


\begin{proof}[Proof of Proposition~\ref{th:relaxed_aadmm_ode}]
The argument is similar to the proof of Proposition~\ref{th:relaxed_admm_ode}.
Combining the optimality condition of \eqref{eq:raadmm2} with \eqref{eq:raadmm3}
we get (note that now $\rho = 1/\epsilon^{2}$)
\begin{equation}\label{eq:uk2}
u_k \in \epsilon^2 \partial g(z_k).
\end{equation}
Combining this and update \eqref{eq:raadmm4} into the optimality condition of \eqref{eq:raadmm1}
yield
\begin{equation} 
\A^T \dfrac{\A x_{k+1} - \hat{z}_k}{\epsilon^2} + \gamma_k \A^T \dfrac{u_k - u_{k-1}}{\epsilon^2}
\in -\partial  f(x_{k+1}) - \A^T \partial g(z_{k}).
\end{equation}
Using update \eqref{eq:raadmm3} we can write this as
\begin{equation} \label{eq:optfirsta}
\alpha^{-1}\A^T \left( 
\dfrac{ u_{k+1} - \hat{u}_k }{\epsilon^2} + 
\dfrac{z_{k+1} - \hat{z}_k}{\epsilon^2}
\right) + \gamma_k \A^T \dfrac{u_k - u_{k-1}}{\epsilon^2}
\in -\partial  f(x_{k+1}) - \A^T \partial g(z_{k}).
\end{equation}

Now we start taking the limit $\epsilon \to 0$, where we recall that $t = \epsilon k$. 
From  the update \eqref{eq:raadmm4} we have
\begin{equation}
\hat{U}(t+\epsilon) = U(t+\epsilon)  + \dfrac{ t + 1}{ t + \epsilon  r + 1} \left( U(t+\epsilon) - U(t) \right),
\end{equation}
which implies   $\hat{U}(t) = U(t)$ in the limit $\epsilon \to 0$.
By the same argument we conclude from \eqref{eq:raadmm5} that 
$\hat{Z}(t) = Z(t)$.
Now from \eqref{eq:uk2} we  have  
$U(t) = \hat{U}(t) = 0$ as $\epsilon \to 0$, and all their derivatives also vanish.
From update \eqref{eq:raadmm5},
\begin{equation}
z_{k+1} - \hat{z}_k = z_{k+1} - z_k - \dfrac{\epsilon (k-1) + 1}{\epsilon (k-1)  + 1 + \epsilon r }
(z_k - z_{k-1}).
\end{equation}
Adding and subtracting $z_{k} - z_{k-1}$ yields
\begin{equation}
\dfrac{z_{k+1} - \hat{z}_k}{\epsilon^2} 
= \dfrac{z_{k+1} - 2 z_k + z_{k-1}}{\epsilon^2} 
+ \dfrac{\epsilon \rr  }{\epsilon (k-1) + 1 + \epsilon r}  \dfrac{ z_k - z_{k-1}}{\epsilon^2} .
\end{equation}
In the limit $\epsilon \to 0$ we thus obtain
\begin{equation}
\dfrac{z_{k+1} - \hat{z}_k}{\epsilon^2} \to \ddot{Z}(t) + \dfrac{r}{t+1} \dot{Z}(t).
\end{equation}
We have a similar result for $(u_{k+1} - \hat{u}_k)/\epsilon^2$, however since
$U(t) = \dot{U}(t) = \ddot{U}(t) = 0$  the first term
in \eqref{eq:optfirsta} vanishes. 
Note also that
\begin{equation}
\gamma_k \dfrac{u_k - u_{k-1} }{\epsilon^2} = \dfrac{\epsilon (k-1) + 1}{\epsilon(k-1) + 1+ \epsilon r} \dfrac{u_k - u_{k-1}}{\epsilon^2}
\end{equation}
so in the limit $\epsilon\to0$ this becomes
\begin{equation}
\gamma_k \dfrac{u_k-u_{k-1}}{\epsilon^2} \to \dfrac{\dot{U}(t)}{\epsilon} \to 0 ,
\end{equation}
since  from \eqref{eq:uk2} we know that $U(t)$ and $\dot{U}(t)$ go to zero with $\epsilon^2$.
Finally, from \eqref{eq:raadmm3} we conclude that $Z(t) = \A X(t)$ in the limit $\epsilon\to 0$.
Therefore, in the limit $\epsilon \to 0$, the discrete inclusion becomes the differential
inclusion
\begin{equation}
\alpha^{-1} \AtA \left( \ddot{X}(t) + \dfrac{r}{t+1} \dot{X}(t) \right) \in
-\partial f(X(t)) - \A^T \partial g (Z(t))  ,
\end{equation}
which along 
with \eqref{lemma} yields \eqref{eq:relaxed_aadmm_ode}.

The first initial condition is obviously $X(0) = x_0$.
For the velocity, consider
the optimality condition for \eqref{eq:raadmm1} with $k=0$, i.e,
$0 \in \partial f(x_{1}) + \rho \A^T(\A x_{1} - \hat{z}_0 + \hat{u}_0)$.
Since $\hat{z}_0 = z_0 = \A x_0$ and $\hat{u}_0 = u_0 = 0$
(recall Remark~\ref{rmk:initial}),
we have $\AtA (x_1 - x_0) \in - \epsilon^2 \partial f(x_1)$.
When $\epsilon \to 0$ this gives $\dot{X}(0) = \tfrac{ X(\epsilon) - X(0) }{ \epsilon} = 0$.
%
\end{proof}

\begin{proof}[Proof of Proposition~\ref{th:heavy_ball_admm_ode}]
The steps are  the same
as in the proof of Proposition~\ref{th:relaxed_aadmm_ode}, 
but with  $\gamma_k = \gamma = 1 - \epsilon \rr$.
\end{proof}

\section{Omitted derivations of convergence rates}
\label{sec:more_rates}

\begin{proof}[Proof of Proposition~\ref{th:relaxed_admm_convex} (convex case)]
Consider
\begin{equation}
\label{eq:lyapunov11}
\E(X, t) \equiv
\alpha t \tilde{\objective}   +
\dfrac{1}{2}\big\| \A \tilde{X} \big\|^2 ,
\end{equation}
where  $X\equiv X_\lambda(t) $ is a trajectory of \eqref{ode1}, 
and $\objective \equiv \objective_\lambda(X_\lambda(t))$; note that we are omitting the
parameter $\lambda$ for the sake of notation.  
Moreover, we use the perturbed variables \eqref{eq:perturb}.
This is a slight modification
of the ``Lyapunov'' function considered in \cite{Franca:2018} and the following steps
are equally similar.
Taking the total time derivative of \eqref{eq:lyapunov11} along trajectories of \eqref{ode1} yields
\begin{equation}\label{Edot-1}
\dot{\E} = - t \big\| \A \dot{\tilde{X}}  \big\|^2
+ \alpha  \big(  \tilde{\objective}  
- \big\langle  \tilde{X}  ,  \nabla \tilde{\objective}  \big\rangle \big) .
\end{equation}
From the convexity of $\objective$, i.e., relation \eqref{eq:strongly_convex} with $\mu=0$,
the second term above is negative. Thus 
$\dot{\E} \le  0$; actually $\dot{\E} < 0$ for $X \ne x^\star$. 
Thus
$\E(X(t), t) \le \E(x_0, 0)$ which implies that
\begin{equation}
\alpha t \big(\objective_\lambda(X_\lambda(t)) - \objective^\star\big)
\le \dfrac{1}{2} \big\| \A (x_0 - x^\star)\big\|^2,
\end{equation}
where we restored $\lambda$. Taking $\lambda \downarrow 0$, we obtain
the  the upper bound \eqref{eq:relaxed_admm_convex} for the associated
differential inclusion \eqref{eq:relaxed_admm_ode}.
\end{proof}

\begin{proof}[Proof of Proposition~\ref{th:relaxed_aadmm_convex} (convex case)]
Consider the regularized differential equation \eqref{ode2}.
We use the centered variables \eqref{eq:perturb} and omit the regularization 
parameter $\lambda$
for simplicity. Thus, consider the function
\begin{equation}
\label{eq:lyapunov3_2}
\E(X, \dot{X}, t) \equiv
\dfrac{\alpha (t+1)^2}{(\rr-1)^2} \tilde{\objective} 
+ \dfrac{1}{2}\left\|\A \left( \tilde{X} +
\dfrac{t+1}{\rr-1}\dot{\tilde{X}} \right)\right\|^2.
\end{equation}
The total time derivative along trajectories of \eqref{ode2} is
\begin{equation} \label{eq:edot-2}
  \dot{\E} = \dfrac{2\alpha(t+1)}{(\rr-1)^2 } \tilde{\objective}
- \dfrac{\alpha (t+1)}{\rr-1} \big\langle \tilde{X}, \nabla \tilde\objective \big\rangle.
\end{equation}
We have
$\big\langle \nabla\tilde \objective,  \tilde{X}  \big\rangle
\ge \tilde{\objective}$
from the convexity of $\objective$, thus
\begin{equation}
\label{eq:aadmm_der_neg}
\dot{\E} \le - \dfrac{\alpha(t+1)(\rr-3)}{(\rr-1)^2  } \tilde{\objective}  \le 0 ,
\end{equation}
where the last inequality follows because $\rr \ge 3$.
Thus  $\E(X(t), \dot{X}(t), t) \le \E(x_0, 0, 0)$, i.e.,
\begin{equation}
\dfrac{\alpha (t+1)^2}{(r-1)^2} \left[ \objective_\lambda(X_\lambda(t)) - \objective^\star \right] 
\le \dfrac{\alpha}{(r-1)^2}\left[ \objective_\lambda(x_0) - x^\star \right] + \dfrac{1}{2} \| \A (x_0 - x^\star)\|^2,
\end{equation}
where we restored the parameter $\lambda$.
This relation  implies \eqref{eq:relaxed_aadmm_convex} in the limit $\lambda \downarrow 0$
with a constant
\begin{equation} \label{constant1}
C \equiv \objective(x_0) - \objective^\star   +
\dfrac{(\rr-1)^2 \sigma_1^2(\A)}{2 \alpha} \| x_0 - x^\star \|^2   .
\end{equation}
\end{proof}

\section{Incorporating nonlinear constraints}
\label{sec:constrained}

It has been stressed in Remarks~\ref{rmk_constr1} and \ref{rmk_constr2} that 
the previous ADMM algorithms are already suitable to handle some types 
of constraints. This is reflected into the modeling differential inclusions
due to the nonsmooth nature of $\partial \phi$ and the preconditioning factor $(\AtA)^{-1}$. 
We now discuss one possible approach
to  extend our methods to account 
for  general \emph{nonlinear equality and inequality constraints}.
Our motivation is the  framework of \cite{Bertsekas:1976}
where constrained optimization problems are solved from solutions to 
unconstrained optimization problems.\footnote{We thank an anonymous 
referee for suggesting this approach.} 
We first discuss how to adapt the previous ADMM variants within this approach,
and then provide a dynamical systems perspective based on singular perturbation
theory.

\subsection{Equality constraints}
Consider the optimization 
problem \eqref{eq:minimize} subject to \emph{nonlinear equality constraints}, 
\begin{equation}\label{eq:minimize_constr}
\begin{split}
&\min_{x} \left\{ \objective(x) \equiv f(x) + g(\A x) \right\} \\
&\mbox{subject to $h_i(x) = 0$, for $i=1,\dotsc,p$} ,
\end{split}
\end{equation}
where $h_i : \mathbb{R}^n \to \mathbb{R}$ is a constraint function.
Denote by $h(x) \equiv (h_1(x), \dotsc, h_p(x))^T$ the vector of constraints, and by $\J_{ij}(x) \equiv \partial h_i(x) / \partial x_j$ its Jacobian
matrix.  Consider the augmented Lagrangian
\begin{equation} \label{eq:aug_lagran}
\mathcal{L}(x, y_\ell, c_\ell) = \objective(x) + \mathcal{L}_{\C}(x, y_\ell, c_\ell) ,
\end{equation}
where
\begin{equation}\label{eq:aug_lagran2}
\mathcal{L}_{\C}(x, y_\ell, c_\ell) \equiv 
\big\langle y_\ell, h(x) \big\rangle  + (c_\ell / 2) \| h(x) \|^2,
\end{equation}
$y_\ell \in \mathbb{R}^p$,
and $\{c_\ell > 0\}$ is a \emph{given sequence} of (increasing or nondecreasing) constants.\footnote{For instance, $c_\ell = s^\ell$ for some $s > 1$. Moreover, one can let $c_\ell = c < \infty$ after a certain number of iterations $\ell \ge \ell^\star$, for sufficiently large  $c$.}  The method consists of sequential unconstrained minimization of 
$\mathcal{L}(x, y_\ell, c_\ell)$. Specifically, for iterations $\ell=0,1,\dotsc$
the general algorithm is given by \cite{Bertsekas:1976}
\begin{subequations}\label{eq:gen_constr}
\begin{align}
x_{\ell+1} &= {\argmin}_x \, \mathcal{L}(x, y_\ell, c_\ell) , \label{eq:constr1} \\
y_{\ell+1} &= y_\ell + c_\ell h(x_{\ell+1}) . \label{eq:constr2}
\end{align} 
\end{subequations}
Formally, $y_\ell$ must be within a bounded
subset $\mathcal{S} \subset \mathbb{R}^p$---$\mathcal{S}$ is arbitrary, but ideally should be as small as possible around the true Lagrange
multiplier $y^\star$. Thus, update \eqref{eq:constr2} holds as long as 
$y_{\ell+1} \in \mathcal{S}$, otherwise $y_{\ell+1} = y_\ell$.

\begin{assumption}\label{th:assump5}
In this section we must rely on the assumptions of \cite[Sec. 2]{Bertsekas:1976}, i.e.,
the critical point $(x^\star, y^\star)$ satisfy standard
2nd order sufficiency conditions for a local minimum of  
problem \eqref{eq:minimize_constr}.  Moreover, the Hessian matrices
$\nabla^2 \objective$, $\nabla^2 h$ are Lipschitz continuous.
\end{assumption}

Thus, in this section we assume sufficient smoothness of all functions
and that the Lagrange multiplier $y^\star$ is unique.
Under these  conditions, 
if the constant $c_\ell \to c < \infty$
then $\lim_{\ell \to \infty} 
\sup \| y_{\ell+1} - y^\star\| / \| y_{\ell} - y^\star \| 
= \bigO(1/c)$, i.e., the convergence is linear 
with a rate proportional to 
$1/c$ \cite[Prop.~1]{Bertsekas:1976} (if $c_k \to \infty$  this rate is superlinear).
This approach therefore allows us to  solve equality constrained problems 
provided the unconstrained problem \eqref{eq:constr1} is solved properly.
One may  use any suitable method for this part,  
such as our previous accelerated variants of ADMM \eqref{eq:relaxed_aadmm}
and \eqref{eq:heavy_ball_admm}. We consider two situations:

\begin{itemize}
\item $g\ne 0$: If the function $g$ is present in  problem \eqref{eq:minimize_constr}
then the part with the constraints must be incorporated
into the proximal operator of $f$, i.e., the only change is in the update
\eqref{eq:raadmm1}/\eqref{eq:heavy_ball_admm1} which now becomes
\begin{equation} \label{eq:gdiffzero}
x_{k+1} = \mbox{$\argmin_x$} 
\big\{ f(x) + \mathcal{L}_C(x, y_\ell, c_\ell)   
+ (\rho/2) \| \A x - \hat{z}_k + u_k \|^2 \big\}.
\end{equation}

\item $g = 0$: If the function $g$ is absent from problem 
\eqref{eq:minimize_constr} then ADMM has a convenient structure that
allows us to solve for the constraints independently of $f$, i.e.,
algorithms \eqref{eq:relaxed_aadmm}
and \eqref{eq:heavy_ball_admm} remain  the same but now with $\A = \bm{I}$, so that
the previous linear constraint of ADMM ensures $z = x$, and we define  the ``new'' $g$ to be $g \equiv \mathcal{L}_C$.
Thus,  update \eqref{eq:raadmm2}/\eqref{eq:heavy_ball_admm2} becomes
\begin{equation}\label{eq:gequalzero}
z_{k+1} = \mbox{$\argmin_z$} \big\{ \mathcal{L}_C(z,  y_\ell, c_\ell)  + (\rho/2) \| \alpha x_{k+1}  + (1-\alpha)\hat{z}_k - z + \hat{u}_k \|^2\big\}.
\end{equation}
\end{itemize}

For a general constraint $h(x)$ the proximal operators
\eqref{eq:gdiffzero} and \eqref{eq:gequalzero} are unlikely to have 
closed form solutions.
Thus, in practice, one may need to devise approximations. 
For instance, since the Jacobian 
$\J_{ij} = \partial h_i / \partial x_j$ has full rank (under Assumption~\ref{th:assump5} the
rows of $\J$ are linear independent) 
we can replace \eqref{eq:gdiffzero} by 
\begin{equation} \label{eq:gdiffzeroap}
x_{k+1} =  \argmin_x \big\{ f(x) + y_\ell^T \J(x_k) + c_\ell h(x_k)^T \J(x_k) 
+ (\rho/2) \| \A x - \hat{z}_k + u_k \|^2  \big\} ,
\end{equation}
and now we are back to the typical case where only the proximal
operator of $f$ is needed. Moreover, even if $h(x)$ is nonconvex,
the proximal operator \eqref{eq:gdiffzeroap} is well-defined for convex $f(x)$.
Similarly,  
we can approximate \eqref{eq:gequalzero}   as
\begin{equation} \label{eq:gequalzeroap}
z_{k+1} = \alpha x_{k+1} + (1-\alpha)\hat{z}_k + \hat{u}_k 
 - (1/\rho)\left(y_\ell^T \J(z_k) + c_\ell h^T(z_k) \J(z_k) \right).
\end{equation} 

\subsection{Inequality constraints}
\emph{Inequality constraints} can be easily transformed into 
equality constraints by  introducing slack variables.
   Consider 
\begin{equation}\label{eq:minimize_constr2}
\begin{split}
&\min_{x} \left\{ \objective(x) \equiv f(x) + g(\A x) \right\} \\
&\mbox{subject to $\underline{h}_i(x) \le 0$, for $i=1,\dotsc,q$},
\end{split}
\end{equation}
where  $\underline{h}_i: \mathbb{R}^n \to \mathbb{R}$, and
$\underline{h}(x) \equiv 
\big(
\underline{h}_1(x), \dotsc, \underline{h}_q(x)\big)^T $. The above problem is equivalent to the equality constrained problem
\begin{equation} \label{eq:minimize_constr22}
\begin{split}
&\min_{x, s} \objective(x)  \\
&\mbox{subject to $\underline{h}_i(x) + s_i = 0$, $s_i \ge 0$},
\end{split}
\end{equation}
for $i=1,\dotsc,q$, and where $s \in \mathbb{R}^q$ is the slack variable.
Thus, setting
$h(x) = \underline{h}(x) + s$ into the Lagrangian
\eqref{eq:aug_lagran}, it is possible to solve for  $s$ explicitly 
\cite[Sec.~5]{Bertsekas:1976}. 
We provide this derivation for convenience.

 Replacing $h(x) = \underline{h}(x) + s$ into
\eqref{eq:aug_lagran}--\eqref{eq:aug_lagran2} yields
\begin{equation}
\mathcal{L}_{\C} = \big\langle y_\ell, \underline{h}(x) + s \big\rangle
+ (c_\ell/2) \big\| \underline{h}(x) + s \big\|^2 .
\end{equation}
Optimizing over $s$ requires
$0 = \nabla_s  \mathcal{L} = y_\ell + c_\ell \big( \underline{h}(x) + s \big)$.
Since we must have $s \ge 0$ the solution is  given by
\begin{equation} 
s = \max\big\{0,  -(1/c_\ell) y_\ell - \underline{h}(x) \big\} ,
\end{equation}
where the $\max$ is taken component-wise, on each entry of the vector.
Update \eqref{eq:constr2} thus becomes
\begin{equation}
y_{\ell+1} = y_\ell + c_\ell (\underline{h}(x) + s) 
= \max\big\{0, y_\ell + c_\ell \underline{h}(x) \big\} .
\end{equation}
When $s = 0$ we have
\begin{equation}
\begin{split}
\mathcal{L}_{\C} &= \big\langle y_\ell, \underline{h}(x) \big\rangle 
+ (c_\ell/2) \big\| \underline{h}(x) \big\|^2  \\
&= \dfrac{1}{2 c_\ell} \left( \big\| y_\ell + c_\ell \underline{h}(x) \big\|^2 - \big\| y_\ell \big\|^2  \right) ,
\end{split}
\end{equation} 
and when $s = -(1/c_\ell) y_\ell - \underline{h}(x)$ we have
\begin{equation}
\mathcal{L}_{\C} = -\dfrac{1}{2 c_\ell} \| y_\ell \|^2 .
\end{equation}
We can combine these last two expressions into
\begin{equation} \label{eq:aug_lagran3}
\mathcal{L}_{\C} = \dfrac{1}{2 c_\ell} \left( \max\left\{0, \big\| y_\ell + c_\ell\underline{h}(x) \big\|^2 \right\} - \big\| y_\ell \big\|^2 \right) .
\end{equation}

Therefore, in the case of inequality constraints \eqref{eq:minimize_constr22}, the method
\eqref{eq:gen_constr} is slightly modified into
\begin{subequations} \label{eq:gen_ineq} 
\begin{align}
x_{\ell + 1} &= {\argmin}_x \mathcal{L}(x, y_\ell, c_\ell), \label{eq:gen_ineq1} \\
y_{\ell+1} &= \max\left\{ 0, y_\ell + c_\ell \underline{h}(x) \right\}. \label{eq:gen_ineq2}
\end{align}
\end{subequations}

With the function \eqref{eq:aug_lagran3}, the modification into the previous accelerated
ADMM methods is the same as for the equality constrained case, i.e.,
\eqref{eq:gdiffzero} and  \eqref{eq:gequalzero}. If $\underline{h}(x)$
is differentiable then one can easily approximate these proximal operators by using
$\nabla \mathcal{L}_{\C}$, in the
same fashion as \eqref{eq:gdiffzeroap} and \eqref{eq:gequalzeroap}. 

Within the above formalism, we 
are now able to  combine both equality and inequality constraints
by applying either \eqref{eq:aug_lagran2}--\eqref{eq:gen_constr}  or 
\eqref{eq:gen_ineq}--\eqref{eq:aug_lagran3} to the respective components
of the problem.  This approach significantly extends the range of applicability
of the previous variants of ADMM, enabling the treatment of
general (smooth) nonlinear constraints.
Finally,  although the framework of \cite{Bertsekas:1976} assumes sufficient smoothness
of the objective function (and constraints), empirically ADMM would operate the same
if the objective function is nonsmooth, however the convergence guarantees of
\cite{Bertsekas:1976} can no longer be granted.

\subsection{Perspective from singular perturbation theory}

We now provide a perspective from
\emph{multiple scale analysis} and \emph{singular perturbation theory} \cite{Verhulst,Kuehn}
to the above constrained optimization framework.
We focus on the case of equality constraints \eqref{eq:gen_constr} since, as discussed,  
inequalities constraints can  be reduced to this case.

We can interpret the method \eqref{eq:gen_constr} as a 
dynamical system evolving in \emph{two different time scales}.  
Because $c_\ell$ is assumed to be large,
the variable $y$ is rapidly-varying, whereas the variable $x$ is slowly-varying---%
note that $c_\ell$ is fixed during the evolution of $x$.
Let us thus introduce a \emph{time-scale parameter}\footnote{Here we use 
$\epsilon$ for the time-scale parameter for consistency with standard notation
in multiple scale analysis ($\epsilon$ is not a discretization step size).}
$\epsilon \sim 1/c_\ell$ so that
$0 < \epsilon  \ll 1$.
First, 
suppose we use the dynamical system \eqref{eq:relaxed_admm_ode} (in the smooth setting)
to describe the minimization problem \eqref{eq:constr1}. A dynamical
model that is somewhat consistent with algorithm \eqref{eq:gen_constr} is  
\begin{subequations} \label{eq:slow_fast1}
\begin{align}
\M \dot{X} &= -\nabla \mathcal{L}(X, Y), \\
\epsilon \dot{Y} &=  h(X) ,
\end{align}
\end{subequations}
where  $\dot{X} \equiv d X / dt$, $\dot{Y} \equiv dY/dt$,
$\M \equiv \alpha^{-1} \AtA$,
\begin{equation} \label{eq:nablaL}
\nabla \mathcal{L} = \nabla \phi(X) +  \J^T(X) Y + c \J^T(X) h(X), 
\end{equation}
and we recall that
$\J_{ij} \equiv \partial h_i / \partial X_j$ the the Jacobian matrix of the constraints.%
\footnote{Note that in \eqref{eq:nablaL} we did not include a parameter $1/\epsilon$
with the last term since it would imply that $X$ is also a fast variable. Instead,
we kept a fixed constant $c > 0$ that can however be large.}
The limit $\epsilon \to 0$ is a \emph{singular limit}, characterized by an
abrupt change of the dynamics.
In this case, the system \eqref{eq:slow_fast1} becomes
the so-called \emph{degenerate system} 
\begin{subequations} \label{eq:fast_subsystem}
\begin{align}
\M \dot{X} &= - \nabla \mathcal{L}(X, Y), \label{fast_sub1} \\ 
0 &= h(X), \label{fast_sub2}
\end{align}
\end{subequations}
which is a \emph{differential algebraic equation}; the evolution occurs in the
\emph{slow manifold} 
$\C \equiv \{ (X, Y) \, | \, h(X) = 0 \}$.
Note that \eqref{fast_sub2} is an equation for $Y$, i.e., its  roots specify
$Y^\star = Y^\star(X)$, and the equation of motion in $\C$ is described by 
$\M \dot{X} = - \nabla \mathcal{L}(X, Y^\star(X))$.
Thus, $(X^\star, Y^\star)$ is a critical point of  the degenerate system
\eqref{eq:fast_subsystem}
if and only if 
\begin{equation} \label{1stopt}
\nabla \objective(X^\star) + \J^T(X^\star) Y^\star = 0, \qquad h(X^\star) = 0,
\end{equation}
which are precisely the 1st order optimality conditions for the
constrained optimization problem.

Now, consider the system \eqref{eq:slow_fast1} evolving on 
the fast time scale $\tau \equiv t /\epsilon$.
We  have 
\begin{subequations}
\begin{align} 
\M X' &= - \epsilon \nabla\mathcal{L}(X, Y),  \\ 
Y' &= h(X),
\end{align}
\end{subequations}
where $X' \equiv dX/d\tau$ and $Y' \equiv dY / d\tau$. In the singular limit  $\epsilon \to 0$ we obtain the so-called \emph{boundary layer system}
\begin{equation} \label{bound_lay}
Y' = h(X)
\end{equation}
with $X' = 0$. Thus, $X$ is constant (i.e., a fixed parameter) on this time scale.
Importantly, note that the update \eqref{eq:constr2} is precisely  a discretization of the
 equation \eqref{bound_lay} (written back to the slow time $t$).
On the other hand, the combination of this step together
with the update \eqref{eq:constr1} emulates the degenerate system \eqref{eq:fast_subsystem},
whose critical points correspond to the 1st order optimality conditions.
This explains why the method \eqref{eq:gen_constr} is expected to yield
a solution  of the constrained optimization problem.
Moreover, it is now clear why the sequence $\{c_\ell > 0\}$ has to increase, i.e., 
it simulates precisely the singular limit $\epsilon \to 0$.

Let us mention that 
there exists a general result in singular perturbation theory due to Tikhonov
(see, e.g., \cite[pp. 433--436]{OMalley:1968}) ensuring that, under appropriate
conditions, such as stability of the boundary layer system with respect to $Y$,%
\footnote{Our case is even simpler since  the vector field $h(X)$ of
\eqref{bound_lay} does not even depend on $Y$.}
solutions of a nonlinear dynamical system
containing a singular parameter $\epsilon$ tend to solutions
of the associated degenerate problem 
when  $\epsilon \to 0$. 

It is clear that we can also replace the other dynamical systems
\eqref{eq:relaxed_aadmm_ode} and \eqref{eq:relaxed_heavy_ball_aadmm_ode}
to model the minimization part encoded by the  update \eqref{eq:constr1}. 
Instead of \eqref{eq:slow_fast1} we now have the full system
\begin{subequations} \label{fast_slow2}
\begin{align}
\dot{V} &= -\gamma(t) V - \M^{-1}\nabla\mathcal{L}(X, Y), \\
\dot{X} &= V, \\
\epsilon \dot{Y} &= h(X),
\end{align}
where we introduced the velocity $V = \dot{X}$ to write the system in 1st order form,
and $\gamma(t) = \rr/(t+1)$ or $\gamma(t) = \rr = \mbox{const.}$.
\end{subequations}
The degenerate system is obtained with $\epsilon \to 0$, i.e.,
\begin{subequations} \label{degen2}
\begin{align}
\M\left[ \ddot{X} + \gamma(t) \dot{X} \right] &= -\nabla \mathcal{L}(X, Y), \\
0 &= h(X),
\end{align}
\end{subequations}
where we wrote back into 2nd order form.
This corresponds to a constrained version of system \eqref{eq:relaxed_heavy_ball_aadmm_ode} (in the smooth case), and the dynamics lie on the slow manifold where the constraints
are satisfied.  In the fast time scale $\tau \equiv t / \epsilon$ the system
\eqref{fast_slow2} tends to the boundary layer system
\begin{equation}
Y' = h(X) ,
\end{equation}
as before, and where now $\dot{V} = \ddot{X} = 0$ and $\dot{X} = 0$, i.e., again
$X$ becomes a constant parameter.
Note that, also in this case,
$(X, V, Y) = (X^\star, 0, Y^\star)$ is a critical point of the degenerate system
\eqref{degen2} if and only if the 1st order optimality conditions
\eqref{1stopt} are satisfied. 

Finally, note that if for some special type of constraints $h(X)$ 
the Lagrangian $\mathcal{L}$ is (strongly) convex in $X$, i.e.,
if one can show that for Lagrange multipliers 
$Y^\star = Y^\star(X)$ solving the constraint condition
$h(X) = 0$ the gradient $\nabla \mathcal{L}(X, Y^\star(X) )$ is (strongly) monotone,
then our previously derived convergence rates for the unconstrained dynamical systems
immediately apply to this setting.

\bibliography{biblio.bib}

\end{document}